\documentclass[final,3p]{elsarticle}
 \usepackage{graphics}
 \usepackage{graphicx}
 \usepackage{epsfig}
\usepackage{amssymb}
 \usepackage{amsthm}
 \usepackage{lineno}
 \usepackage{amsmath}
   \numberwithin{equation}{section}
\usepackage{mathrsfs}

\newtheorem{thm}{Theorem}[section]
\newtheorem{cor}[thm]{Corollary}
\newtheorem{lem}[thm]{Lemma}
\newtheorem{prop}[thm]{Proposition}
\newtheorem{defn}[thm]{Definition}

\biboptions{sort&compress,square}
\begin{document}
\begin{frontmatter}
\author[rvt1]{ Kai Hua Bao\corref{cor1}}
\cortext[cor1]{Corresponding author.}
\ead{baokh505@nenu.edu.cn}
\author[rvt2]{kun ming hu}
%\author[rvt3]{Jian Wang}
\address[rvt1]{School of Mathematics, Ineer Mongolia University for Nationnalities, TongLiao, 028005, P.R.China}
\address[rvt2]{School of Mathematics and Statistics, Northeast Normal University,
Changchun, 130024, P.R.China}
\title{A Kastler-Kalau-Walze Type Theorem for $7$-dimensional \\ Manifolds with  Boundary about Witten deformation}
\begin{abstract}
In this paper, we give a brute-force proof of the  Kastler-Kalau-Walze type theorem
for $7$-dimensional manifolds with  boundary about Witten deformation, and give a theoritic explaination of the gravitational action for $7$ dimensional manifolds with boundary.
\end{abstract}
\end{frontmatter}
\section{Introduction}
\label{1}
The noncommutative residue plays a prominent role in noncommutative geometry \cite{Gu}\cite{Wo}. Connes \cite{Co1} used the noncommutative
 residue to derive a conformal 4-dimensional Polyakov action analogy. Connes \cite{Co2} proved that the noncommutative
residue on a compact manifold $M$ coincided with the Dixmier's trace on pseudodifferential operators of order $-{\rm {dim}}M$. Several
years ago, Connes made a challenging observation that the noncommutative residue of the square of the inverse of the Dirac
operator was proportional to the Einstein-Hilbert action, which we call the Kastler-Kalau-Walze theorem. Kastler\cite{Ka} gave a
brute-force proof of this theorem. Kalau and Walze \cite{KW} proved this theorem in the normal coordinates system simultaneously.
Ackermann \cite{Ac} gave a note on a new proof of this theorem by means of the heat kernel expansion.

  Recently, Ponge defined lower dimensional volumes of Riemannian manifolds by the Wodzicki residue \cite{Po}.
Fedosov et al defined a noncommutative residue on Boutet de Monvel's algebra and proved that it was a
unique continuous trace \cite{FGLS}. Wang generalized  the Connes' results to the case of manifolds with boundary in \cite{Wa1} \cite{Wa2} ,
and proved a Kastler-Kalau-Walze type theorem for the Dirac operator and the signature operator for 3, 4-
dimensional manifolds with boundary\cite{Wa3}. Wang also generalized the definition of lower dimensional volumes to manifolds with boundary, and found a Kastler-Kalau-Walze type theorem for higher dimensional manifolds with boundary\cite{Wa4}. Weiping Zhang introduced an elliptic differential operator-Witten deformation in \cite{wpz}. In\cite{KH1,KH2}, we proved Kastler-Kalau-Walze type theorem for Witten deformation for 4, 6-
dimensional manifolds with boundary\cite{Wa3}. Furthermore, we consider higher dimensional case. The motivation of this paper is to establish a Kastler-Kalau-Walze type theorem associated with Witten deformation for $7$-dimensional manifolds with boundary, and give a theoritic explaination of the gravitational action for  $7$-dimensional manifolds with boundary.

This paper is organized as follows: In Section 2, we define lower dimensional volumes of compact Riemannian manifolds
with  boundary. In Section 3, for $7$-dimensional spin compact manifolds with boundary and the associated Witten deformation,
we compute lower dimensional volumes $ \widetilde{{\rm Wres}}[(\pi^+D_{T}^{-2})^2]$ and get a Kastler-Kalau-Walze type theorem in this case. We also give a theoritic explaination of the gravitational action for $7$ dimensional manifolds with boundary.

\section{Lower-Dimensional volumes of compact manifolds with boundary about Witten deformation }
 In this section we consider an $n$-dimensional oriented compact Riemannian manifold $(M, g^{M})$ with boundary $\partial_{ M}$ equipped
with a fixed spin structure. We assume that the metric $g^{M}$ on $M$ has
the following form near the boundary
 \begin{equation}
 g^{M}=\frac{1}{h(x_{n})}g^{\partial M}+\texttt{d}x _{n}^{2} ,
\end{equation}
where $g^{\partial M}$ is the metric on $\partial M$. Let $U\subset
M$ be a collar neighborhood of $\partial M$ which is diffeomorphic $\partial M\times [0,1)$. By the definition of $h(x_n)\in C^{\infty}([0,1))$
and $h(x_n)>0$, there exists $\tilde{h}\in C^{\infty}((-\varepsilon,1))$ such that $\tilde{h}|_{[0,1)}=h$ and $\tilde{h}>0$ for some
sufficiently small $\varepsilon>0$. Then there exists a metric $\hat{g}$ on $\hat{M}=M\bigcup_{\partial M}\partial M\times
(-\varepsilon,0]$ which has the form on $U\bigcup_{\partial M}\partial M\times (-\varepsilon,0 ]$
 \begin{equation}
\hat{g}=\frac{1}{\tilde{h}(x_{n})}g^{\partial M}+\texttt{d}x _{n}^{2} ,
\end{equation}
such that $\hat{g}|_{M}=g$.
We fix a metric $\hat{g}$ on the $\hat{M}$ such that $\hat{g}|_{M}=g$.

Firstly, we will recall the expression of Witten deformation $D_{T}$ and $D_{T}^{2}$ near the boundary\cite{KH1,KH2}. Let $\nabla^L$ denote the Levi-Civita connection about $g^M$. In the local coordinates $\{x_i; 1\leq i\leq n\}$ and the
 fixed orthonormal frame $\{\widetilde{e_1},\cdots,\widetilde{e_n}\}$, the connection matrix $(\omega_{s,t})$ is defined by
\begin{equation}
\nabla^L(\widetilde{e_1},\cdots,\widetilde{e_n})= (\widetilde{e_1},\cdots,\widetilde{e_n})(\omega_{s,t}).
\end{equation}
 Let $\epsilon (\widetilde{e_j*} ),~\iota (\widetilde{e_j*} )$ be the exterior and interior multiplications respectively. Write
\begin{equation}
c(\widetilde{e_j})=\epsilon (\widetilde{e_j*} )-\iota (\widetilde{e_j*} );~~
\bar{c}(\widetilde{e_j})=\epsilon (\widetilde{e_j*} )+\iota
(\widetilde{e_j*} ).
\end{equation}
 The Witten deformation is defined by
\begin{equation}
D_{T}=d+\delta+T\bar{c}(V)=\sum^n_{i=1}c(\widetilde{e_i})[\widetilde{e_i}+\frac{1}{4}\sum_{s,t}\omega_{s,t}
(\widetilde{e_i})[\bar{c}(\widetilde{e_s})\bar{c}(\widetilde{e_t})-c(\widetilde{e_s})c(\widetilde{e_t})]]+T\bar{c}(V),
\end{equation}
where $d,\delta$,$V\in\Gamma(TM)$,  any $T\in{\bf R}$.

By proposition 4.6 of \cite{wpz}, we have
\begin{equation}
D^{2}_{T}=(d+\delta)^{2}+\sum_{i}^{n}c(\widetilde{e_i})\nabla_{\widetilde{e_i}}^{TM}V+T^{2}|V|^{2}.
\end{equation}

By \cite{Y}, $(d+\delta)^{2}$ is expressed by
\begin{equation}
(d+\delta)^{2}
=-\triangle_{0}-\frac{1}{8}\sum_{ijkl}R_{ijkl}\bar{c}(\widetilde{e_i})\bar{c}(\widetilde{e_j})c(\widetilde{e_k})c(\widetilde{e_l})-\frac{1}{4}s.
\end{equation}

Let $g^{ij}=g(dx_{i},dx_{j})$, $\xi=\sum_{k}\xi_{j}dx_{j}$ and $\nabla^L_{\partial_{i}}\partial_{j}=\sum_{k}\Gamma_{ij}^{k}\partial_{k}$, we denote
\begin{eqnarray}
\sigma_{i}=-\frac{1}{4}\sum_{s,t}\omega_{s,t}
(\widetilde{e_i})c(\widetilde{e_s})c(\widetilde{e_t})
;~~~a_{i}=\frac{1}{4}\sum_{s,t}\omega_{s,t}
(\widetilde{e_i})\bar{c}(\widetilde{e_s})\bar{c}(\widetilde{e_t}).
\end{eqnarray}

Denote
\begin{eqnarray}
\xi^{j}=g^{ij}\xi_{i},~~~~\Gamma^{k}=g^{ij}\Gamma_{ij}^{k},~~~~\sigma^{j}=g^{ij}\sigma_{i},~~~~a^{j}=g^{ij}a_{i},
\end{eqnarray}
then $D_{T}$ can be written as
\begin{equation}
D_{T}=\sum^n_{i=1}c(\widetilde{e_i})(\widetilde{e_i}+\sigma_{i}+a_{i})+T\bar{c}(V).
\end{equation}

By \cite{Ac}, \cite{Y}, we have
\begin{equation}
-\triangle_{0}=\Delta=-g^{ij}(\nabla^L_{i}\nabla^L_{j}-\Gamma_{ij}^{k}\nabla^L_{k}),
\end{equation}
then by (2.7) we have
\begin{eqnarray}
D^{2}_{T}&=&-\sum_{i,j}g^{i,j}\Big[\partial_{i}\partial_{j}+2\sigma_{i}\partial_{j}+2a_{i}\partial_{j}-\Gamma_{i,j}^{k}\partial_{k}+(\partial_{i}\sigma_{j})+
(\partial_{i}a_{j})+\sigma_{i}\sigma_{j}+\sigma_{i}a_{j}+a_{i}\sigma_{j}+a_{i}a_{j} -\Gamma_{i,j}^{k}\sigma_{k}\nonumber\\
&&-\Gamma_{i,j}^{k}a_{k}\Big]
-\frac{1}{8}\sum_{ijkl}R_{ijkl}\bar{c}(\widetilde{e_i})\bar{c}(\widetilde{e_j})c(\widetilde{e_k})c(\widetilde{e_l})+\frac{1}{4}s+\sum_{i}^{n}c(\widetilde{e_i})\bar{c}(\nabla_{\widetilde{e_i}}^{TM}V)+T^{2}|V|^{2}.
\end{eqnarray}

To define the lower dimensional volume, some basic facts and formulae about Boutet de Monvel's calculus which can be found in Sec.2 in \cite{Wa3}
are needed.
 Let  \begin{equation*}
  F:L^2({\bf R}_t)\rightarrow L^2({\bf R}_v);~F(u)(v)=\int e^{-ivt}u(t)dt
  \end{equation*}
   denote the Fourier transformation and
$\Phi(\overline{{\bf R}^+}) =r^+\Phi({\bf R})$ (similarly define $\Phi(\overline{{\bf R}^-}$)), where $\Phi({\bf R})$
denotes the Schwartz space and
  \begin{equation*}
r^{+}:C^\infty ({\bf R})\rightarrow C^\infty (\overline{{\bf R}^+});~ f\rightarrow f|\overline{{\bf R}^+};~
 \overline{{\bf R}^+}=\{x\geq0;x\in {\bf R}\}.
\end{equation*}
We define $H^+=F(\Phi(\overline{{\bf R}^+}));~ H^-_0=F(\Phi(\overline{{\bf R}^-}))$ which are orthogonal to each other. We have the following
 property: $h\in H^+~(H^-_0)$ iff $h\in C^\infty({\bf R})$ which has an analytic extension to the lower (upper) complex
half-plane $\{{\rm Im}\xi<0\}~(\{{\rm Im}\xi>0\})$ such that for all nonnegative integer $l$,
 \begin{equation*}
\frac{d^{l}h}{d\xi^l}(\xi)\sim\sum^{\infty}_{k=1}\frac{d^l}{d\xi^l}(\frac{c_k}{\xi^k}),
\end{equation*}
as $|\xi|\rightarrow +\infty,{\rm Im}\xi\leq0~({\rm Im}\xi\geq0)$.

 Let $H'$ be the space of all polynomials and $H^-=H^-_0\bigoplus H';~H=H^+\bigoplus H^-.$ Denote by $\pi^+~(\pi^-)$ respectively the
 projection on $H^+~(H^-)$. For calculations, we take $H=\widetilde H=\{$rational functions having no poles on the real axis$\}$ ($\tilde{H}$
 is a dense set in the topology of $H$). Then on $\tilde{H}$,
 \begin{equation}
\pi^+h(\xi_0)=\frac{1}{2\pi i}\lim_{u\rightarrow 0^{-}}\int_{\Gamma^+}\frac{h(\xi)}{\xi_0+iu-\xi}d\xi,
\end{equation}
where $\Gamma^+$ is a Jordan close curve
included ${\rm Im}(\xi)>0$ surrounding all the singularities of $h$ in the upper half-plane and
$\xi_0\in {\bf R}$. Similarly, define $\pi^{'}$ on $\tilde{H}$,
\begin{equation}
\pi'h=\frac{1}{2\pi}\int_{\Gamma^+}h(\xi)d\xi.
\end{equation}
So, $\pi'(H^-)=0$. For $h\in H\bigcap L^1(R)$, $\pi'h=\frac{1}{2\pi}\int_{R}h(v)dv$ and for $h\in H^+\bigcap L^1(R)$, $\pi'h=0$.

Denote by $\mathcal{B}$ Boutet de Monvel's algebra, we recall the main theorem in \cite{FGLS}.
\begin{thm}\label{th:32}{\bf(Fedosov-Golse-Leichtnam-Schrohe)}
 Let $X$ and $\partial X$ be connected, ${\rm dim}X=n\geq3$,
 $A=\left(\begin{array}{lcr}\pi^+P+G &   K \\
T &  S    \end{array}\right)$ $\in \mathcal{B}$ , and denote by $p$,$s$ and $b$ the local symbols of pseudo-differential operators $P$ and $S$, Singular Green operator $G$ respectively, $T$ is trace operator.
 Define:
 \begin{eqnarray}
{\rm{\widetilde{Wres}}}(A)&=&\int_X\int_{\bf S}{\rm{tr}}_E\left[p_{-n}(x,\xi)\right]\sigma(\xi)dx \nonumber\\
&&+2\pi\int_ {\partial X}\int_{\bf S'}\left\{{\rm tr}_E\left[({\rm{tr}}b_{-n})(x',\xi')\right]+{\rm{tr}}
_F\left[s_{1-n}(x',\xi')\right]\right\}\sigma(\xi')dx',
\end{eqnarray}
Then~~ a) ${\rm \widetilde{Wres}}([A,B])=0 $, for any
$A,B\in\mathcal{B}$;~~ b) It is a unique continuous trace on
$\mathcal{B}/\mathcal{B}^{-\infty}$.
\end{thm}
 Let $p_{1},p_{2}$ be nonnegative integers and $p_{1}+p_{2}\leq n$. From Sec.2.1 of \cite{Wa3}, we have the following definition.
\begin{defn}{\rm (\cite{WJ1})} Lower-dimensional volumes of compact spin manifolds with boundary about Witten deformation $D_{T}$ is defined by
   \begin{equation}
  {\rm Vol}^{(p_1,p_2)}_nM:=\widetilde{{\rm Wres}}[\pi^+D_{T}^{-p_1}\circ\pi^+D_{T}^{-p_2}].
\end{equation}
\end{defn}
Denote by $\sigma_{l}(A)$ the $l$-order symbol of an operator A. An application of (2.1.4) in \cite{Wa1} shows that
\begin{equation}
\widetilde{{\rm Wres}}[\pi^+D_{T}^{-p_1}\circ\pi^+D_{T}^{-p_2}]=\int_M\int_{|\xi|=1}{\rm
trace}_{S(TM)}[\sigma_{-n}(D_{T}^{-p_1}\circ D_{T}^{-p_2})]\sigma(\xi)\texttt{d}x+\int_{\partial
M}\Phi,
\end{equation}
\begin{eqnarray}
\Phi&=&\int_{|\xi'|=1}\int^{+\infty}_{-\infty}\sum^{\infty}_{j, k=0}
\sum\frac{(-i)^{|\alpha|+j+k+1}}{\alpha!(j+k+1)!}
 {\rm trace}_{S(TM)}
\Big[\partial^j_{x_n}\partial^\alpha_{\xi'}\partial^k_{\xi_n}
\sigma^+_{r}D_{T}^{-p_1}(x',0,\xi',\xi_n)\nonumber\\
&&\times\partial^\alpha_{x'}\partial^{j+1}_{\xi_n}\partial^k_{x_n}\sigma_{l}
D_{T}^{-p_2}(x',0,\xi',\xi_n)\Big]d\xi_n\sigma(\xi')\texttt{d}x',
\end{eqnarray}
and the sum is taken over $r-k+|\alpha|+\ell-j-1=-n,r\leq-p_{1},\ell\leq-p_{2}$.
%In this paper will compute
%\begin{eqnarray}
%\widetilde{{\rm Wres}}\big[\pi^+(D_{T}^{*})^{-1}\circ \pi^+D_{T}^{-(2\overline{n}-1)}\big],~~~~~~~~for~~n=(2\overline{n}+2);\nonumber\\
%\widetilde{{\rm Wres}}\big[\pi^+(D_{T}^{*})^{-2}\circ \pi^+D_{T}^{-(2\overline{n})}\big],~~~~~~~~~~~~for~~n=(2\overline{n}+4).\nonumber
%\end{eqnarray}

The following proposition is the key of the computation of lower-dimensional volumes of compact spin
manifolds with boundary.
\begin{prop}{\rm(\cite{Wa4})}  The following identity holds:
 \begin{eqnarray}
&& 1)~  When~ p_1+p_2=n,~ then, ~ {\rm Vol}^{(p_1,p_2)}_nM=c_0{\rm Vol}_M;\\
&& 2)~  when~  p_1+p_2\equiv n ~ {\rm mod}~  1, ~{\rm Vol}^{(p_1,p_2)}_nM=\int_{\partial M}\Phi.
\end{eqnarray}
\end{prop}

\section{ A Kastler-Kalau-Walze type theorem for $7$-dimensional spin manifolds with boundary }
In this section, we compute the lower dimensional volume ${\rm Vol}^{(2,2)}_7$ for $7$-dimensional spin compact manifolds with boundary and prove a
Kastler-Kalau-Walze type theorem in this case. By Proposition 2.4, we have
 \begin{equation}
 \widetilde{{\rm Wres}}[\pi^+D_{T}^{-2}\circ\pi^+D_{T}^{-2}]=\int_{\partial M}\Phi.
\end{equation}
So we only need to compute $\int_{\partial M}\Phi$.

Firstly, we will compute some symbols of Witten deformation. Let
\begin{eqnarray}
\Gamma^{k}=\sum_{i,j<n}\sum_{l<n}g^{ij}g^{lk}\langle \nabla^{ L}_{\partial_{i}}\partial_{j}, \partial_{l}\rangle
 +\sum_{l<n}g^{lk}\langle \nabla^{ L}_{\partial_{n}}\partial_{n}, \partial_{l}\rangle.
\end{eqnarray}
By conclusion of \cite{KH1,KH2}, we have
\begin{lem}\label{le:31}{\rm(\cite{KH1,KH2})}
The symbols of the Witten deformation are
\begin{eqnarray}
\sigma_{-1}(D^{-1}_{T})&=&\frac{\sqrt{-1}c(\xi)}{|\xi|^2}; \nonumber\\
\sigma_{-2}(D^{-1}_{T})&=&\frac{c(\xi)\sigma_{0}(D_{T})c(\xi)}{|\xi|^4}+\frac{c(\xi)}{|\xi|^6}\sum_jc(dx_j)
\Big[\partial_{x_j}(c(\xi))|\xi|^2-c(\xi)\partial_{x_j}(|\xi|^2)\Big] ;\nonumber\\
\sigma_{-2}(D^{-2}_{T})&=&|\xi|^{-2}; \nonumber\\
\sigma_{-3}(D^{-2}_{T})&=&-\sqrt{-1}|\xi|^{-4}\xi_{k}(\Gamma^{k}-2a^{k}-2\sigma^{k})-\sqrt{-1}|\xi|^{-6}2\xi^j\xi_\alpha\xi_\beta\partial_jg^{\alpha\beta}.
\end{eqnarray}
\end{lem}
Now we will compute symbol $\sigma_{-4}(D_{T}^{-2})$.
Since the equation {\rm (4.16)} of {\rm \cite{KW}} is
\begin{eqnarray}
\sigma^{\tilde{\triangle}^{-1}}_{-4}(x,\xi)=\sigma^{-1}_{2}(\gamma^{2}_{-1}+\gamma_{-2})+i\sigma^{-2}_{2}\partial_{\xi_{\mu}}\sigma_{2}\partial_{x^{\mu}}\gamma_{-1},
\end{eqnarray}
where \begin{eqnarray}
\sigma^{\tilde{\triangle}^{1}}(x,\xi)&=&\sigma_{2}+\sigma_{1}+\sigma_{0},\nonumber\\
\gamma_{-1}(x,\xi)&=&-\sigma^{-1}_{2}\sigma_{1}-i\sigma^{-2}_{2}\partial_{\xi_{\mu}}\sigma_{2}\partial_{x^{\mu}}\sigma_{2},\nonumber\\
\gamma_{-2}(x,\xi)&=&-\sigma^{-1}_{2}\sigma_{0}-\sigma^{-2}_{2}(i\partial_{\xi_{\mu}}\sigma_{1}\partial_{x^{\mu}}\sigma_{2}+\frac{1}{2}\partial_{\xi_{\mu}}\partial_{\xi_{\nu}}\sigma_{2}\partial_{x^{\mu}}\partial_{x^{\nu}}\sigma_{2})\nonumber\\
&~&+\sigma^{-3}_{2}\partial_{\xi_{\mu}}\partial_{\xi_{\nu}}\sigma_{2}\partial_{x^{\mu}}\sigma_{2}\partial_{x^{\nu}}\sigma_{2}.\nonumber\\
\end{eqnarray}

Then by Lemma 3.1 and some calculations we get
\begin{lem}\label{le:32}
\begin{eqnarray}
&&\sigma_{-4}(D_{T}^{-2})=\sigma_{-4}(D^{-2})-4i|\xi|^{-6}(\Gamma^{k}-2\sigma^{k})\xi_{k}a^{l}\xi_{l}-4|\xi|^{-6}(a^{k}\xi_{k})^{2}+4|\xi|^{-8}a^{k}\xi_{k}\partial_{\xi_{\mu}}(|\xi|^{2})\partial_{x^{\mu}}(|\xi|^{2})\nonumber\\
&&+|\xi|^{-4}(\partial^{i}a_{j}+\sigma^{i}a_{j}+a^{i}\sigma_{j}+a^{i}a_{j}-\Gamma^{k}a_{k}
+\frac{1}{8}\sum_{ijkl}R_{ijkl}\bar{c}(\widetilde{e_i})\bar{c}(\widetilde{e_j})c(\widetilde{e_k})c(\widetilde{e_l})-\sum_{i}^{n}c(\widetilde{e_i})\bar{c}(\nabla_{\widetilde{e_i}}^{TM}V) \nonumber\\
&&-T^{2}|V|^{2})-2|\xi|^{-6}\partial_{\xi_{\mu}}(a^{k}\xi_{k})\partial_{x^{\mu}}(|\xi|^{2})-2|\xi|^{-2}\partial_{\xi_{\mu}}(|\xi|^{2})\partial_{x^{\mu}}(|\xi|^{2})a^{k}\xi_{k}
-2|\xi|^{-4}\partial_{\xi_{\mu}}(|\xi|^{2})\partial_{x^{\mu}}(a^{k}\xi_{k}),\nonumber\\
\end{eqnarray}
Where  $\sigma_{-4}(D^{-2})$ has the following expression by the equation {\rm (115) }in {\rm\cite{WJ2}}
\begin{eqnarray}
\sigma_{-4}(D^{-2})&=&-|\xi|^{-6}\xi_{k}\xi_{l}(\Gamma^{k}-2\sigma^{k})(\Gamma^{l}-2\sigma^{l})
   +2|\xi|^{-8}\xi^{k}\xi_{l}\xi_{\alpha}\xi_{\beta}(\Gamma^{l}-2\sigma^{l})\partial^{x}_{\mu}g^{\alpha\beta}\nonumber\\
   &&+|\xi|^{-4}(\partial^{x_{k}}\sigma_{k}+\sigma^{k}\sigma_{k}-\Gamma^{k}\sigma_{k} )-\frac{1}{4}|\xi|^{-4}s(x)
  -2|\xi|^{-6}\xi^{k}\xi_{l}\partial^{x}_{k}(\Gamma^{l}-2\sigma^{l})\nonumber\\
   && +12|\xi|^{-10}\xi^{k}\xi_{l}\xi_{\alpha}\xi_{\beta}\xi_{\gamma}\xi_{\delta}\partial^{x}_{k}g^{\alpha\beta}
   \partial^{x}_{l}g^{\gamma\delta}
  -4|\xi|^{-8}\xi^{k}\xi_{\alpha}\xi_{\gamma}\xi_{\delta}\partial^{x}_{k}g^{l\alpha}\partial^{x}_{l}g^{\gamma\delta}\nonumber\\
  && -4|\xi|^{-8}\xi^{k}\xi^{l}\xi_{\gamma}\xi_{\delta}\partial^{x}_{kl}g^{\gamma\delta}
  +|\xi|^{-6}\xi_{\alpha}\xi_{\beta}(\Gamma^{k}-2\sigma^{k})\partial^{x}_{k}g^{\alpha\beta}
   -|\xi|^{-6}\xi_{\alpha}\xi_{\beta}g^{kl}\partial^{x}_{kl}g^{\alpha\beta}\nonumber\\
   &&+2|\xi|^{-8}\xi_{\alpha}\xi_{\beta}\xi_{\gamma}\xi_{\delta}g^{kl}\partial^{x}_{k}g^{\alpha\beta}\partial^{x}_{l}g^{\gamma\delta}.
\end{eqnarray}
\end{lem}

Since $\Phi$ is a global form on $\partial M$, so for any fixed point $x_{0}\in\partial M$, we can choose the normal coordinates
$U$ of $x_{0}$ in $\partial M$(not in $M$) and compute $\Phi(x_{0})$ in the coordinates $\widetilde{U}=U\times [0,1)$ and the metric
$\frac{1}{h(x_{n})}g^{\partial M}+\texttt{d}x _{n}^{2}$. The dual metric of $g^{M}$ on $\widetilde{U}$ is
$h(x_{n})g^{\partial M}+\texttt{d}x _{n}^{2}.$ Write
$g_{ij}^{M}=g^{M}(\frac{\partial}{\partial x_{i}},\frac{\partial}{\partial x_{j}})$;
$g^{ij}_{M}=g^{M}(d x_{i},dx_{j})$, then

\begin{equation}
[g_{i,j}^{M}]=
\begin{bmatrix}\frac{1}{h( x_{n})}[g_{i,j}^{\partial M}]&0\\0&1\end{bmatrix};\quad
[g^{i,j}_{M}]=\begin{bmatrix} h( x_{n})[g^{i,j}_{\partial M}]&0\\0&1\end{bmatrix},
\end{equation}
and
\begin{equation}
\partial_{x_{s}} g_{ij}^{\partial M}(x_{0})=0,\quad 1\leq i,j\leq n-1;\quad g_{i,j}^{M}(x_{0})=\delta_{ij}.
\end{equation}

Let $\{E_{1},\cdots, E_{n-1}\}$ be an orthonormal frame field in $U$ about $g^{\partial M}$ which is parallel along geodesics and
$E_{i}=\frac{\partial}{\partial x_{i}}(x_{0})$, then $\{\widetilde{E_{1}}=\sqrt{h(x_{n})}E_{1}, \cdots,
\widetilde{E_{n-1}}=\sqrt{h(x_{n})}E_{n-1},\widetilde{E_{n}}=dx_{n}\}$ is the orthonormal frame field in $\widetilde{U}$ about $g^{M}.$
Locally $S(TM)|\widetilde{U}\cong \widetilde{U}\times\wedge^{*}_{C}(\frac{n}{2}).$ Let $\{f_{1},\cdots,f_{n}\}$ be the orthonormal basis of
$\wedge^{*}_{C}(\frac{n}{2})$. Take a spin frame field $\sigma: \widetilde{U}\rightarrow Spin(M)$ such that
$\pi\sigma=\{\widetilde{E_{1}},\cdots, \widetilde{E_{n}}\}$ where $\pi: Spin(M)\rightarrow O(M)$ is a double covering, then
$\{[\sigma, f_{i}], 1\leq i\leq 6\}$ is an orthonormal frame of $S(TM)|_{\widetilde{U}}.$ In the following, since the global form $\Phi$
is independent of the choice of the local frame, so we can compute $\texttt{tr}_{S(TM)}$ in the frame $\{[\sigma, f_{i}], 1\leq i\leq 6\}$.
Let $\{\hat{E}_{1},\cdots,\hat{E}_{n}\}$ be the canonical basis of $R^{n}$ and
$c(\hat{E}_{i})\in cl_{C}(n)\cong Hom(\wedge^{*}_{C}(\frac{n}{2}),\wedge^{*}_{C}(\frac{n}{2}))$ be the Clifford action. By \cite{Wa3}, then

\begin{equation}
c(\widetilde{E_{i}})=[(\sigma,c(\hat{E}_{i}))]; \quad c(\widetilde{E_{i}})[(\sigma, f_{i})]=[\sigma,(c(\hat{E}_{i}))f_{i}]; \quad
\frac{\partial}{\partial x_{i}}=[(\sigma,\frac{\partial}{\partial x_{i}})],
\end{equation}
then we have $\frac{\partial}{\partial x_{i}}c(\widetilde{E_{i}})=0$ in the above frame.

Nextly, we will give some conclusions as our computing tools. By \cite{WJ2} we have
\begin{lem}\label{le:32}
With the metric $g^{M}$ on $M$ near the boundary
\begin{eqnarray}
\partial_{x_j}(|\xi|_{g^M}^2)(x_0)&=&\left\{
       \begin{array}{c}
        0,  ~~~~~~~~~~ ~~~~~~~~~~ ~~~~~~~~~~~~~{\rm if }~j<n; \\[2pt]
       h'(0)|\xi'|^{2}_{g^{\partial M}},~~~~~~~~~~~~~~~~~~~~{\rm if }~j=n.
       \end{array}
    \right. \\
\partial_{x_j}(c(\xi))(x_0)&=&\left\{
       \begin{array}{c}
      0,  ~~~~~~~~~~ ~~~~~~~~~~ ~~~~~~~~~~~~~{\rm if }~j<n;\\[2pt]
\partial_{x_{n}} (c(\xi'))(x_{0}), ~~~~~~~~~~~~~~~~~{\rm if }~j=n,
       \end{array}
    \right.
\end{eqnarray}
where $\xi=\xi'+\xi_{n}\texttt{d}x_{n}$.
\end{lem}

\begin{lem}\label{le:32}{\rm(\cite{WJ2})}
With the metric $g^{M}$ on $M$ near the boundary
\begin{eqnarray}
\partial_{x_i}\partial_{x_j}(|\xi|_{g^M}^2)(x_0)\Big|_{|\xi'|=1}&=&\left\{
       \begin{array}{c}
        0,~~~ ~~~~~~~~~~ ~~~~~~~~~~ ~~~~~~~~~~~~~{\rm if }~i<n,j=n;or~i=n,j<n ;\\
        -\frac{1}{3}\sum_{\alpha ,\beta <n}\Big(R^{\partial_{M}}_{i\alpha j\beta}(x_0)+R^{\partial_{M}}_{i\beta j\alpha}(x_0)\Big)\xi_{\alpha}\xi_{\beta},~~~~~~~~~~~{\rm if }~i,j<n;\\
       h''(0),~~~~~~~~~~~~~~~~~~~~~~~~~~~~~~~~~~~~~~~~~~~~~~~~~~~~~{\rm if }~i=j=n,
       \end{array}
    \right. \\
\partial_{x_i}\partial_{x_j}[c(\xi)](x_0)\Big|_{|\xi'|=1}&=&\left\{
       \begin{array}{c}
       0,~~~ ~~~~~~~~~~ ~~~~~~~~~~ ~~~~~~~~~~~~~{\rm if }~i<n,j=n;or~i=n,j<n ;\\
        \frac{1}{6}\sum_{l ,t <n}\xi_{l}\Big(R^{\partial_{M}}_{i\alpha j\beta}(x_0)+R^{\partial_{M}}_{i\beta j\alpha}(x_0)\Big)c(\widetilde{e_{t}}),~~~~~~~~~~~~~{\rm if }~i,j<n;\\
       \bigg(\frac{3}{4}(h'(0))^{2}-\frac{1}{2}h''(0)\bigg)\sum_{j<n}\xi_{j}c(\widetilde{e_{j}}),~~~~~~~~~~~~~~~~~~~~~{\rm if }~j=n,
       \end{array}
    \right.
\end{eqnarray}
where $\xi=\xi'+\xi_{n}\texttt{d}x_{n}$.
\end{lem}
Similar to concliusions in\cite{WJ2}, we get
\begin{lem}\label{le:32}
With the metric $g^{M}$ on $M$ near the boundary
\begin{eqnarray}
\Gamma^{k}(x_0)&=&\left\{
       \begin{array}{c}
        0  ~~~~~~~~~~ ~~~~~~~~~~ ~~~~~~~~~~~~ ~~~~~~~~~~~~ ~~~~~~~~~~~~{\rm if }~k<n; \\[2pt]
       3h'(0)~~~~~~~~~~~~~~~~~~~~~~~~~~ ~~~~~~~~~~~~ ~~~~~~~~~~~~{\rm if }~k=n.
       \end{array}
    \right. \\
\sigma^{k}(x_0)&=&\left\{
       \begin{array}{c}
     -\frac{1}{4}h'(0)c(\widetilde{e_{k}})c(\widetilde{e_{n}})  ~~~~~~~~~~~~~~~~ ~~~~~~~~~~~~~~~~~~~{\rm if }~k<n;\\[2pt]
0 ~~~~~~~~~~~~~~~~~~~~~ ~~~~~~~~~~~~ ~~~~~~~~~~~~~~~ ~~~~~~~{\rm if }~k=n.
       \end{array}
    \right.\\
  a^{k}(x_0)&=&\left\{
       \begin{array}{c}
     \frac{1}{4}h'(0)\overline{c}(\widetilde{e_{k}})\overline{c}(\widetilde{e_{n}})  ~~~~~~~~~~~~~~~~ ~~~~~~~~~~~~~~~~~~~~~{\rm if }~k<n;\\[2pt]
0 ~~~~~~~~~~~~~~~~~~~~~ ~~~~~~~~~~~~ ~~~~~~~~~~~~~~~ ~~~~~~~~{\rm if }~k=n.
       \end{array}
    \right.\\
\end{eqnarray}
\end{lem}
Similar to concliusions in\cite{WJ2}, we get
\begin{lem}\label{le:32}
With the metric $g^{M}$ on $M$ near the boundary
\begin{eqnarray*}
\partial_{x_\gamma}\Gamma^{k}(x_0)\Big|_{|\xi'|=1}&=&\left\{
       \begin{array}{c}
        \frac{5}{6} \sum_{i <n}R^{\partial_{M}}_{i\gamma ik}(x_0),
        ~~~~~~~~~~~~~~~~~~~~~~~~~  ~~~~~~~~~~~~~~{\rm if }~ \gamma<n,k<n; \\  \\
          ~~~~0,~~~~~~~~~~~~~~~~~~~~~~~~~~~~~~~~~~~~~~~~~~~~~~~~~
           ~~~~~{\rm if }~\gamma<n,k=n;\\[2pt]
    \\
        ~~~~ 0,~~~~~~~~~~~~~~~~~~~~~~~~~~~~~~~~~~~~~~~~~~~~~~~~~
           ~~~~~~{\rm if }~\gamma=n,k<n;\\[2pt]   \\
3h''(0)-\frac{9}{2}\big(h'(0)\big)^{2}, ~~~~~~~~~~~~~~~~~~~~~~~~~~~~~~~~~~~~~~{\rm if }~\gamma=n,k=n.
       \end{array}
    \right. \\  \nonumber\\
\partial_{x_\gamma}\sigma^{k}(x_0)\Big|_{|\xi'|=1}&=&\left\{
       \begin{array}{c}
    -\frac{1}{8} \sum_{s\neq t <n}R^{\partial_{M}}_{k\gamma st}(x_0)c(\widetilde{e_{s}})c(\widetilde{e_{t}}),
        ~~~~~~~~~~~~~~~~~~~~~~~{\rm if }~ \gamma<n,k<n; \\  \\
          ~~~~0,~~~~~~~~~~~~~~~~~~~~~~~~~~~~~~~~~~~~~~~~~~~~~~~~~
           ~~~~~~{\rm if }~\gamma<n,k=n;\\[2pt]
    \\
       -\sum_{t <n}\Big(\frac{3}{8}\big(h'(0)\big)^{2}-\frac{1}{4}h''(0)\Big)c(\widetilde{e_{n}})c(\widetilde{e_{t}}),
        ~~~~~~~~~~~~~{\rm if }~\gamma=n,k<n;\\[2pt]   \\
 -\frac{1}{8}\sum_{t <n}\Big(\big(h'(0)\big)^{2}-h''(0)\Big)c(\widetilde{e_{s}})c(\widetilde{e_{t}}).
  ~~~~~~~~~~~~~~~~{\rm if }~\gamma=n,k=n.
       \end{array}
    \right. \\  \nonumber\\
\partial_{x_\gamma}a^{k}(x_0)\Big|_{|\xi'|=1}&=&\left\{
       \begin{array}{c}
    \frac{1}{8} \sum_{s\neq t <n}R^{\partial_{M}}_{k\gamma st}(x_0)\overline{c}(\widetilde{e_{s}})\overline{c}(\widetilde{e_{t}}),
        ~~~~~~~~~~~~~~~~~~~~~~~~~~{\rm if }~ \gamma<n,k<n; \\  \\
          ~~~~0,~~~~~~~~~~~~~~~~~~~~~~~~~~~~~~~~~~~~~~~~~~~~~~~~~
           ~~~~~~~~{\rm if }~\gamma<n,k=n;\\[2pt]
    \\
       \sum_{t <n}\Big(\frac{3}{8}\big(h'(0)\big)^{2}-\frac{1}{4}h''(0)\Big)\overline{c}(\widetilde{e_{n}})\overline{c}(\widetilde{e_{t}}),
        ~~~~~~~~~~~~~~~{\rm if }~\gamma=n,k<n;\\[2pt]   \\
 \frac{1}{8}\sum_{t <n}\Big(\big(h'(0)\big)^{2}-h''(0)\Big)\overline{c}(\widetilde{e_{s}})\overline{c}(\widetilde{e_{t}}).
  ~~~~~~~~~~~~~~~~~~{\rm if }~\gamma=n,k=n.
       \end{array}
    \right. \\\nonumber\\
\end{eqnarray*}
\end{lem}

Now we will compute $\Phi$ (see formula (2.15) for definition of $\Phi$).
Since the sum is taken over $-r-\ell+1+k+j+|\alpha|=7, \ r, \ell\leq-2$, then we have the $\int_{\partial_{M}}\Phi$ is the sum of the following
 fifteen cases:\\

\textbf{Case (1)}: \ $r=-2, \ \ell=-2, \ k=0, \ j=1, \ |\alpha|=1$\\

From (2.15), we have
\begin{equation}
\text{ Case \ (1)}=\frac{i}{2}\int_{|\xi'|=1}\int_{-\infty}^{+\infty}\sum_{|\alpha|=1}\text{trace}
\Big[\partial_{x_{n}}\partial_{\xi'}^{\alpha}\pi_{\xi_{n}}^{+}\sigma_{-2}(D_{T}^{-2})
\partial_{x'}^{\alpha}\partial_{\xi_{n}}^{2}\sigma_{-2}(D_{T}^{-2})\Big](x_{0})\texttt{d}\xi_{n}\sigma(\xi')\texttt{d}x' .
\end{equation}
By Lemma 3.1 and Lemma 3.3, for $i<n$, we have
\begin{equation}
\partial_{x_i}\sigma_{-2}(D_{T}^{-2})(x_0)=\partial_{x_i}\big(|\xi|^{-2}\big)(x_0)=0.
\end{equation}
So Case (1) vanishes.

\textbf{Case (2)}: \  $r=-2, \ \ell=-2, \ k=0, \ j=2, \ |\alpha|=0$\\

From (2.15), we have
\begin{equation}
\text{ Case \ (2)}=\frac{i}{6}\int_{|\xi'|=1}\int_{-\infty}^{+\infty}\sum_{j=2}\text{trace}\Big[\partial_{x_{n}}^{2}\pi_{\xi_{n}}^{+}
\sigma_{-2}(D_{T}^{-2})\partial_{\xi_{n}}^{3}\sigma_{-2}(D_{T}^{-2})\Big](x_{0})\texttt{d}\xi_{n}\sigma(\xi')\texttt{d}x'.
\end{equation}
By Lemma 3.1 and Lemma 3.3 and a simple calculation, we get
\begin{equation}
\partial_{\xi_{n}}^{3}\sigma_{-2}(D_{T}^{-2})(x_{0})\Big|_{|\xi'|=1}
=\frac{24\xi_{n}-24\xi_{n}^{3}}{(1+\xi_{n}^{2})^{4}},
\end{equation}
and
\begin{equation}
\partial^{2}_{x_{n}}\sigma_{-2}(D_{T}^{-2})(x_{0})
=\frac{2(h'(0))^{2}}{(1+\xi_{n}^{2})^{3}}-\frac{h''(0)}{(1+\xi_{n}^{2})^{2}}.
\end{equation}
Since $\partial^{2}_{x_{n}}$ and $\pi^+_{\xi_n}$ can be exchange, and by (2.13) we have
\begin{equation}
\partial_{x_{n}}^{2}\pi_{\xi_{n}}^{+}\sigma_{-2}(D_{T}^{-2})(x_{0})\Big|_{|\xi'|=1}
=\frac{-3i\xi_n^{2}-9\xi_n+8i }{8(\xi_n-i)^3}(h'(0))^{2}+\frac{2+i\xi_n }{4(\xi_n-i)^2}h''(0).
\end{equation}
Note that, for $7$-dimensional spin compact manifolds with boundary we have $\texttt{trace}_{S(TM)}[\texttt{id}]=8$. Then by (3.22),  (3.24) and some direct computations, we obtain
\begin{eqnarray}
&&\text{trace}\Big[\partial_{x_{n}}^{2}\pi_{\xi_{n}}^{+}
\sigma_{-2}(D_{T}^{-2})\partial_{\xi_{n}}^{3}\sigma_{-2}(D_{T}^{-2})\Big](x_{0})\Big|_{|\xi'|=1}\nonumber\\
&&=\big(h'(0)\big)^{2}\frac{(-3i\xi_{n}^{2}-9\xi_{n}+8i)(24\xi_{n}-24\xi_{n}^{3})}{(\xi_{n}-i)^{3}(1+\xi_{n}^{2})^{4}}
 +h''(0)\frac{(4+2i\xi_{n})(24\xi_{n}-24\xi_{n}^{3})} {(\xi_{n}-i)^{2}(1+\xi_{n}^{2})^{4}}.
\end{eqnarray}
Therefore
\begin{eqnarray}
\text{ Case \ (2) }
&=&\frac{i}{6}\big(h'(0)\big)^{2}  \frac{ 2\pi i }{6!}
\bigg[\frac{72i\xi_{n}^{5}+216\xi_{n}^{4}-264i\xi_{n}^{3}-216\xi_{n}^{2}+192i\xi_{n}}{(\xi_{n}+i)^{4}}\bigg]^{(6)}\bigg|_{\xi_{n}=i}\Omega_{5}\texttt{d}x'  \nonumber\\
&&+\frac{i}{6}h''(0) \frac{ 2\pi i }{5!}
\bigg[\frac{-48i\xi_{n}^{4}-96\xi_{n}^{3}+48i\xi_{n}^{2}+96\xi_{n}} {(\xi_{n}^{2}+i)^{4}}\bigg]^{(5)}\bigg|_{\xi_{n}=i}\Omega_{5}\texttt{d}x'  \nonumber\\
&=&\Big(\frac{7}{8}\big(h'(0)\big)^{2}-\frac{3}{8}h''(0)\Big)\pi\Omega_{5}\texttt{d}x',
\end{eqnarray}
where $\Omega_{5}$ is the canonical volume of $S^{5}$.

\textbf{Case (3)}: \ $r=-2, \ \ell=-2, \ k=0, \ j=0, \ |\alpha|=2$

From (2.15), we have
\begin{equation}
\text{ Case \ (3)}=\frac{i}{2}\int_{|\xi'|=1}\int_{-\infty}^{+\infty}\sum_{|\alpha|=2}
\text{trace}\Big[\partial_{\xi'}^{\alpha}\pi_{\xi_{n}}^{+}\sigma_{-2}(D_{T}^{-2})
\partial_{x'}^{\alpha}\partial_{\xi_{n}}\sigma_{-2}(D_{T}^{-2})\Big](x_{0})\texttt{d}\xi_{n}\sigma(\xi')\texttt{d}x' .
\end{equation}
By Lemma 3.5 and a simple calculation we have
\begin{eqnarray}
\partial_{\xi'}^{\alpha}\sigma_{-2}(D_{T}^{-2})(x_{0})\Big|_{|\xi'|=1}
&=&\sum_{i,j<n}\partial_{\xi_{j}}\partial_{\xi_{i}}\sigma_{-2}(D_{T}^{-2})(x_{0})\Big|_{|\xi'|=1}  \nonumber\\
&=&\sum_{i,j<n}\frac{-2\delta_{i}^{j}}{(1+\xi_{n}^{2})^{2}}+\sum_{i,j<n}\frac{8}{(1+\xi_{n}^{2})^{3}}\xi_{i}\xi_{j}.
\end{eqnarray}
By (3.28)and (2.13) we obtain
\begin{equation}
\pi_{\xi_{n}}^{+}\partial_{\xi'}^{\alpha}\sigma_{-2}(D_{T}^{-2})(x_{0})\Big|_{|\xi'|=1}
=\sum_{i,j<n}\frac{(2+i\xi_{n})\delta_{i}^{j}}{2(\xi_{n}-i)^{2}}+\sum_{i,j<n}\frac{-3i\xi_{n}^{2}-9\xi_{n}+8i}{2(\xi_{n}-i)^{3}}\xi_{i}\xi_{j}.
\end{equation}
On the other hand, by Lemma 3.2 and Lemma 3.5, we obtain
\begin{eqnarray}
\partial_{x'}^{\alpha}\sigma_{-2}(D_{T}^{-2})(x_{0})\Big|_{|\xi'|=1}
&=&\sum_{i,j<n}\frac{-1}{(1+\xi_{n}^{2})^{2}}\partial_{x_{i}}\partial_{x_{j}}(|\xi|^{2})(x_{0})
+\sum_{i,j<n}\frac{2}{(1+\xi_{n}^{2})^{3}}\partial_{x_{j}}(|\xi|^{2})\partial_{x_{i}}(|\xi|^{2})(x_{0}) \nonumber\\
&=&\frac{1}{3(1+\xi_{n}^{2})^{2}}
 \sum_{i,j,\alpha ,\beta <n}\Big(R^{\partial_{M}}_{i\alpha j\beta}(x_0)+R^{\partial_{M}}_{i\beta j\alpha}(x_0)\Big)\xi_{\alpha}\xi_{\beta}
+\frac{2\big(h'(0)\big)^{2}}{(1+\xi_{n}^{2})^{3}}.
\end{eqnarray}
Hence in this case,
\begin{equation}
\partial_{x'}^{\alpha}\partial_{\xi_{n}}\sigma_{-2}(D_{T}^{-2})(x_{0})\Big|_{|\xi'|=1}
=\frac{-4\xi_{n}}{3(1+\xi_{n}^{2})^{3}}
 \sum_{i,j,\alpha ,\beta <n}\Big(R^{\partial_{M}}_{i\alpha j\beta}(x_0)+R^{\partial_{M}}_{i\beta j\alpha}(x_0)\Big)\xi_{\alpha}\xi_{\beta}
+\frac{-12\xi_{n}\big(h'(0)\big)^{2}}{(1+\xi_{n}^{2})^{4}}.
\end{equation}
By (3.29), (3.31) and some direct computations, we obtain
\begin{eqnarray}
&&\text{trace}\Big[\partial_{\xi'}^{\alpha}\pi_{\xi_{n}}^{+}\sigma_{-2}(D_{T}^{-2})
\partial_{x'}^{\alpha}\partial_{\xi_{n}}\sigma_{-2}(D_{T}^{-2})\Big](x_{0})\nonumber\\
&&=\frac{-4\xi_{n}-2i\xi_{n}^{2}}{3(\xi_{n}-i)^{2}(1+\xi_{n}^{2})^{3}}
 \sum_{i,j,\alpha ,\beta <n}\Big(R^{\partial_{M}}_{i\alpha j\beta}(x_0)+R^{\partial_{M}}_{i\beta j\alpha}(x_0)\Big)\xi_{\alpha}\xi_{\beta}\nonumber\\
&&+\frac{-16i\xi_{n}+18\xi_{n}^{2}+6i\xi_{n}^{3}}{3(\xi_{n}-i)^{3}(1+\xi_{n}^{2})^{3}}
 \sum_{i,j,\alpha ,\beta <n}\Big(R^{\partial_{M}}_{i\alpha j\beta}(x_0)+R^{\partial_{M}}_{i\beta j\alpha}(x_0)\Big)
 \xi_{i}\xi_{j}\xi_{\alpha}\xi_{\beta}\nonumber\\
&&+\big(h'(0)\big)^{2}\frac{-12\xi_{n}-6i\xi_{n}^{2}}{(\xi_{n}-i)^{2}(1+\xi_{n}^{2})^{4}}
+\big(h'(0)\big)^{2}\frac{-48i\xi_{n}+54\xi_{n}^{2}+18i\xi_{n}^{3}}{(\xi_{n}-i)^{3}(1+\xi_{n}^{2})^{4}} .
\end{eqnarray}
Similar to (16) in \cite{Ka}, we have
\begin{equation}
\int \xi^{\mu}\xi^{\nu}=\frac{1}{6}[^{\mu\nu}],~~\int \xi^{\mu}\xi^{\nu}\xi^{\alpha}\xi^{\beta}=c_{0}[^{\mu\nu\alpha\beta}],
\end{equation}
where $[^{\mu\nu\alpha\beta}]$ stands for the sum of products of $g^{\alpha\beta}$ determined by all "pairings"
of $\mu\nu\alpha\beta$ and $c_{0}$ is a constant.
Using the integration over $S^{5}$ and the shorthand $ \int=\frac{1}{\pi^{3}}\int_{S^{5}}d^{5}\nu$, we obtain  $\Omega_{5}=\pi^{3} $.
Let $s_{\partial_{M}}$ is the scalar curvature of $\partial_{M}$, then
\begin{equation}
\sum_{i,\alpha, j,\beta <n}R^{\partial_{M}}_{i\alpha j \beta }(x_0)\int_{|\xi'|=1} \xi_{\alpha}\xi_{\beta}\xi_{i}\xi_{j} \sigma(\xi')
 =c\pi^{3}\sum_{i,\alpha, j,\beta <n}R^{\partial_{M}}_{i\alpha j \beta }(x_0)
\Big(\delta_{\alpha}^{\beta}\delta_{i}^{j}+\delta_{\alpha}^{i}\delta_{\beta}^{j}
 +\delta_{\alpha}^{j}\delta_{\beta}^{i}\Big)=0,
\end{equation}
where $c$ is a constant.
Therefore
\begin{eqnarray}
\text{ Case \ (3) }
&=&\frac{i}{2}\Omega_{5}\Big(s_{\partial_{M}}
\int_{-\infty}^{+\infty}\frac{-4\xi_{n}-2i\xi_{n}^{2}}{9(\xi_{n}-i)^{2}(1+\xi_{n}^{2})^{3}}
\texttt{d}\xi_{n}
+\big(h'(0)\big)^{2}
\int_{-\infty}^{+\infty}\frac{4i\xi_{n}-9\xi_{n}^{2}-3i\xi_{n}^{3}}{(\xi_{n}-i)^{3}(1+\xi_{n}^{2})^{4}}\texttt{d}\xi_{n}\Big)\texttt{d}x'\nonumber\\
&=&\frac{i}{2}\Omega_{5}\Big(s_{\partial_{M}}
\frac{ 2\pi i }{4!}\bigg[\frac{-4\xi_{n}-2i\xi_{n}^{2}}
{9(\xi_{n}+i)^{3}}\bigg]^{(4)}\bigg|_{\xi_{n}=i}
+\big(h'(0)\big)^{2}
\frac{ 2\pi i }{6!}\bigg[\frac{4i\xi_{n}-9\xi_{n}^{2}-3i\xi_{n}^{3}}
{(\xi_{n}+i)^{4}}\bigg]^{(6)}\bigg|_{\xi_{n}=i}\Big)\texttt{d}x' \nonumber\\
&=&\frac{\pi}{6}s_{\partial_{M}}\Omega_{5}\texttt{d}x'+\frac{11\pi}{128}\big(h'(0)\big)^{2}\Omega_{5}\texttt{d}x',
\end{eqnarray}
where $\sum_{t,l <n}R^{\partial_{M}}_{tltl}(x_0) $ is the scalar curvature $s_{\partial_{M}}$.

\textbf{Case (4)}: \ $r=-2, \ \ell=-2, \ k=1, \ j=1, \ |\alpha|=0$

From (2.15) and the Leibniz rule, we obtain
\begin{equation}
\text{ Case \ (4)}=\frac{i}{6}\int_{|\xi'|=1}\int_{-\infty}^{+\infty}
\text{trace}\Big[\partial_{x_{n}}\partial_{\xi_{n}}\pi_{\xi_{n}}^{+}\sigma_{-2}(D_{T}^{-2})
\partial_{\xi_{n}}^{2}\partial_{x_{n}}\sigma_{-2}(D_{T}^{-2})\Big](x_{0})\texttt{d}\xi_{n}\sigma(\xi')\texttt{d}x' .
\end{equation}
By Lemma 3.1 and Lemma 3.3, some calculations we obtain
\begin{equation}
\partial_{x_{n}}\partial_{\xi_n}\pi^+_{\xi_n}\sigma_{-2}(D_{T}^{-2})(x_0)|_{|\xi'|=1}=h'(0)\frac{-3-i\xi_n}{4(\xi_n-i)^{3}},
\end{equation}

\begin{equation}
\partial_{\xi_{n}}^{2}\partial_{x_{n}}\sigma_{-2}(D_{T}^{-2})(x_0)|_{|\xi'|=1}
=h'(0)\frac{4-20\xi_{n}^{2}}{(1+\xi_{n}^{2})^{4}}.
\end{equation}
Note that $\texttt{trace}[\texttt{id}]=8$, then by (3.37),  (3.38) and some direct computations, we obtain
\begin{eqnarray}
\text{ Case \ (4) }
&=&\frac{i}{6}\big(h'(0)\big)^{2}\int_{|\xi'|=1}\int_{-\infty}^{+\infty}
\frac{-24-8i\xi_{n}+120\xi_{n}^{2}+40i\xi_{n}^{3}} {(\xi_{n}-i)^{3}(1+\xi_{n}^{2})^{4}}
\texttt{d}\xi_{n}\sigma(\xi')\texttt{d}x'\nonumber\\
&=&\frac{i}{6}\big(h'(0)\big)^{2}  \frac{ 2\pi i }{6!}\bigg[
\frac{-24-8i\xi_{n}+120\xi_{n}^{2}+40i\xi_{n}^{3}}
{(\xi_{n}+i)^{4}}\bigg]^{(6)}\bigg|_{\xi_{n}=i}\Omega_{5}\texttt{d}x'  \nonumber\\
&=&-\frac{5}{8}\big(h'(0)\big)^{2}\pi\Omega_{5}\texttt{d}x'.
\end{eqnarray}

\textbf{Case (5)}: \ $r=-2, \ \ell=-2, \ k=1, \ j=0, \ |\alpha|=1$

From (2.15), we have
\begin{equation}
\text{ Case \ (5)}=\frac{i}{2}\int_{|\xi'|=1}\int_{-\infty}^{+\infty}\sum_{|\alpha|=1}
\text{trace}\Big[\partial_{\xi'}^{\alpha}\partial_{\xi_{n}}\pi_{\xi_{n}}^{+}\sigma_{-2}(D_{T}^{-2})
\partial_{x'}^{\alpha}\partial_{\xi_{n}}\partial_{x_{n}}\sigma_{-2}(D_{T}^{-2})\Big](x_{0})\texttt{d}\xi_{n}\sigma(\xi')\texttt{d}x' .
\end{equation}
By Lemma 3.1 and Lemma 3.3, for $i<n$, we obtain
\begin{eqnarray}
\partial_{x'}\partial_{x_{n}}\sigma_{-2}(D_{T}^{-2})(x_{0})\Big|_{|\xi'|=1}
&=&\frac{-1}{(1+\xi_{n}^{2})^{2}}\partial_{x_{i}}\partial_{x_{n}}(|\xi|^{2})(x_{0})
+\frac{2}{(1+\xi_{n}^{2})^{3}}\partial_{x_{n}}(|\xi|^{2})\partial_{x_{i}}(|\xi|^{2})(x_{0}) \nonumber\\
&=&0.
\end{eqnarray}
Therefore Case (5) vanishes.

\textbf{Case (6)}: \ $r=-2, \ \ell=-2, \ k=2, \ j=0, \ |\alpha|=0$

From (2.15), we have
\begin{equation}
\text{ Case \ (6)}=\frac{i}{6}\int_{|\xi'|=1}\int_{-\infty}^{+\infty}\sum_{k=2}
\text{trace}\Big[\partial_{\xi_{n}}^{2}\pi_{\xi_{n}}^{+}\sigma_{-2}(D_{T}^{-2})
\partial_{\xi_{n}}\partial_{x_{n}}^{2}\sigma_{-2}(D_{T}^{-2})\Big](x_{0})\texttt{d}\xi_{n}\sigma(\xi')\texttt{d}x'.
\end{equation}
By Lemma 3.1, Lemma 3.3, and some calculations, we have
\begin{equation}
\partial_{\xi_{n}}^{2}\pi_{\xi_{n}}^{+}\sigma_{-2}(D_{T}^{-2})(x_{0})\Big|_{|\xi'|=1}
=\frac{-i}{(\xi_{n}-i)^{3}},
\end{equation}
\begin{equation}
\partial_{\xi_{n}}\partial_{x_{n}}^{2}\sigma_{-2}(D_{T}^{-2})(x_{0})\Big|_{|\xi'|=1}
= \frac{4\xi_{n}h''(x_{0})}{(1+\xi_{n}^{2})^{3}}+\frac{-12\xi_{n}(h'(0))^{2}}{(1+\xi_{n}^{2})^{4}}.
\end{equation}
Note that $\texttt{trace}[\texttt{id}]=8$, then by (3.43),  (3.44) and some direct computations, we obtain
\begin{eqnarray}
\text{ Case \ (6) }
&=&\frac{i}{6}h''(0)\int_{|\xi'|=1}\int_{-\infty}^{+\infty}
\frac{-32i\xi_{n}} {(\xi_{n}-i)^{3}(1+\xi_{n}^{2})^{3}}\texttt{d}\xi_{n}\sigma(\xi')\texttt{d}x'\nonumber\\
&&+\frac{i}{6}\big(h'(0)\big)^{2}\int_{|\xi'|=1}\int_{-\infty}^{+\infty}\frac{96i\xi_{n}}
 {(\xi_{n}-i)^{3}(1+\xi_{n}^{2})^{4}}\texttt{d}\xi_{n}\sigma(\xi')\texttt{d}x'\nonumber\\
&=&\frac{i}{6}h''(0)  \frac{ 2\pi i }{5!}\bigg[\frac{-32i\xi_{n}}{(\xi_{n}+i)^{3}}\bigg]^{(5)}\bigg|_{\xi_{n}=i}\Omega_{5}\texttt{d}x'
+\frac{i}{6}\big(h'(0)\big)^{2} \frac{ 2\pi i }{6!}\bigg[\frac{96i\xi_{n}}
{(\xi_{n}+i)^{4}}\bigg]^{(6)}\bigg|_{\xi_{n}=i}\Omega_{5}\texttt{d}x'  \nonumber\\
&=&\Big(-\frac{3}{8}h''(0)+\frac{7}{8}\big(h'(0)\big)^{2}\Big)\pi\Omega_{5}\texttt{d}x'.
\end{eqnarray}

\textbf{Case (7)}: \ $r=-2, \ \ell=-3, \ k=0, \ j=1, \ |\alpha|=0$

From (2.15) and the Leibniz rule, we obtain
\begin{eqnarray}
\text{ Case \ (7)}&=&\frac{1}{2}\int_{|\xi'|=1}\int_{-\infty}^{+\infty}
\text{trace}\Big[\partial_{\xi_{n}}\partial_{x_{n}}\pi_{\xi_{n}}^{+}\sigma_{-2}(D_{T}^{-2})
\partial_{\xi_{n}}\sigma_{-3}(D_{T}^{-2})\Big](x_{0})\texttt{d}\xi_{n}\sigma(\xi')\texttt{d}x' \nonumber\\
  &=&-\frac{1}{2}\int_{|\xi'|=1}\int_{-\infty}^{+\infty}
\text{trace}\Big[\partial_{\xi_{n}}^{2}\partial_{x_{n}}\pi_{\xi_{n}}^{+}\sigma_{-2}(D_{T}^{-2})
\sigma_{-3}(D_{T}^{-2})\Big](x_{0})\texttt{d}\xi_{n}\sigma(\xi')\texttt{d}x'.
\end{eqnarray}
By Lemma 3.1, Lemma 3.3 and some calculations, we have
\begin{equation}
\pi^+_{\xi_n}\partial_{x_n}\sigma_{-2}(D_{T}^{-2})(x_0)|_{|\xi'|=1}=h'(0)\frac{2+i\xi_n}{4(\xi_n-i)^{2}},
\end{equation}
Then
\begin{equation}
\partial_{\xi_{n}}^{2}\pi^+_{\xi_n}\partial_{x_n}\sigma_{-2}(D_{T}^{-2})(x_0)|_{|\xi'|=1}=h'(0)\frac{4+i\xi_n}{2(\xi_n-i)^{4}}.
\end{equation}
In the normal coordinate we have
\begin{eqnarray}
g^{ij}(x_0)&=&\delta_i^j, \nonumber\\
\partial_{x_j}(g^{\alpha\beta})(x_0)&=&\left\{
       \begin{array}{c}
        0  ~~~~~~~~~~ ~~~~ {\rm if }~j<n; \\[2pt]
      h'(0)\delta^\alpha_\beta~~~~~~{\rm if }~j=n.
       \end{array}
    \right.
\end{eqnarray}
So by some calculations we have
\begin{eqnarray}
\sigma_{-3}(D_{T}^{-2})(x_0)|_{|\xi'|=1}
&=&-i|\xi|^{-4}\xi_k(\Gamma^k-2\delta^k)(x_0)|_{|\xi'|=1}-i|\xi|^{-6}2\xi^j\xi_\alpha\xi_\beta
\partial_jg^{\alpha\beta}(x_0)|_{|\xi'|=1}-2i|\xi|^{-4}a^k\xi_k\nonumber\\
&=&\frac{ih'(0)\sum_{k<n}\xi_k
c(\widetilde{e_k})c(\widetilde{e_n})}{2(1+\xi_n^2)^2}+\frac{i3h'(0)\xi_n}{(1+\xi_n^2)^2}\nonumber\\
&-&\frac{2ih'(0)\xi_n}{(1+\xi_n^2)^3}-\frac{ih'(0)\sum_{k<n}\xi_k\overline{c}(\widetilde{e_k})\overline{c}(\widetilde{e_n})}{2(1+\xi_n^2)^2}.\nonumber\\
\end{eqnarray}
We note that
\begin{equation}
\int_{|\xi'|=1}\xi_1\cdots\xi_{2q+1}\sigma(\xi')=0.
\end{equation}

So the first term and the last term in (3.50) has no contribution for computing case
(7), which we will omit in following equation. Combining (3.48), (3.50) and some direct computations, we obtain
\begin{eqnarray}
\text{trace}\Big[\partial_{\xi_{n}}^{2}\partial_{x_{n}}\pi_{\xi_{n}}^{+}\sigma_{-2}(D_{T}^{-2})\sigma_{-3}(D_{T}^{-2})\Big](x_{0})
&=&(h'(0))^{2}\frac{-80i\xi_{n}+20\xi_{n}^{2}-48i\xi_{n}^{3}+12\xi_{n}^{4} } {(\xi_{n}-i)^{4}(1+\xi_{n}^{2})^{3}}.
\end{eqnarray}
Note that $\texttt{trace}[\texttt{id}]=8$, by some direct computations, we obtain
\begin{eqnarray}
\text{ Case \ (7) }
&=&-\frac{1}{2}\big(h'(0)\big)^{2}\int_{|\xi'|=1}\int_{-\infty}^{+\infty}
\big[\frac{-80i\xi_{n}+20\xi_{n}^{2}-48i\xi_{n}^{3}+12\xi_{n}^{4} } {(\xi_{n}-i)^{4}(1+\xi_{n}^{2})^{3}}\nonumber\\
&=&-\frac{1}{2}\big(h'(0)\big)^{2}  \frac{ 2\pi i }{6!}\bigg[
\frac{-80i\xi_{n}+20\xi_{n}^{2}-48i\xi_{n}^{3}+12\xi_{n}^{4}}
{(\xi_{n}+i)^{3}}\bigg]^{(6)}\bigg|_{\xi_{n}=i}\Omega_{5}\texttt{d}x'  \nonumber\\
&=&\frac{21}{8}\big(h'(0)\big)^{2}\pi\Omega_{5}\texttt{d}x'.
\end{eqnarray}

\textbf{Case (8)}: \ $r=-2, \ \ell=-3, \ k=0, \ j=0, \ |\alpha|=1$

From (2.15) and the Leibniz rule, we obtain
\begin{eqnarray}
\text{ Case \ (8)}&=&-\int_{|\xi'|=1}\int_{-\infty}^{+\infty}\sum_{|\alpha|=1}
\text{trace}\Big[\partial_{\xi'}^{\alpha}\pi_{\xi_{n}}^{+}\sigma_{-2}(D_{T}^{-2})
\partial_{x'}^{\alpha}\partial_{\xi_{n}}\sigma_{-3}(D_{T}^{-2})\Big](x_{0})\texttt{d}\xi_{n}\sigma(\xi')\texttt{d}x' \nonumber\\
  &=&\int_{|\xi'|=1}\int_{-\infty}^{+\infty}\sum_{|\alpha|=1}
\text{trace}\Big[\partial_{\xi_{n}}\partial_{\xi'}^{\alpha}\pi_{\xi_{n}}^{+}\sigma_{-2}(D_{T}^{-2})
\partial_{x'}^{\alpha}\sigma_{-3}(D_{T}^{-2})\Big](x_{0})\texttt{d}\xi_{n}\sigma(\xi')\texttt{d}x'.
\end{eqnarray}
By Lemma 3.1, Lemma 3.3 and some calculations, we get
\begin{equation}
\partial_{\xi'}^{\alpha}\sigma_{-2}(D_{T}^{-2})(x_{0})\Big|_{|\xi'|=1}
=\sum_{k<n}\partial_{\xi_{k}}\sigma_{-2}(D_{T}^{-2})(x_{0})\Big|_{|\xi'|=1} =\sum_{k<n}\frac{-2\xi_{k} }{(1+\xi_{n}^{2})^{2}}.
\end{equation}
By some calculations, we obtain
\begin{equation}
\partial_{\xi_{n}}\partial_{\xi'}^{\alpha}\pi_{\xi_{n}}^{+}\sigma_{-2}(D_{T}^{-2})(x_{0})\Big|_{|\xi'|=1}
 =\sum_{k<n}\frac{-3-i\xi_{n} }{2(\xi_{n}-i)^{3}}\xi_{k}.
\end{equation}
By Lemma 3.1, Lemma 3.6 and some direct computations, we obtain
\begin{eqnarray}
\partial_{x'}\sigma_{-3}(D_{T}^{-2})(x_{0})\Big|_{|\xi'|=1}&=&\frac{-5i\xi_{k} }{6(1+\xi_{n}^{2})^{2}}\sum_{ i<n}R^{\partial_{M}}_{ i \gamma i k}(x_0) \nonumber\\
&&+ \frac{i\xi_{k} }{4(1+\xi_{n}^{2})^{2}}\sum_{ s\neq t<n}R^{\partial_{M}}_{ k \gamma st}(x_0)
c(\widetilde{e_{s}})c(\widetilde{e_{t}}) \nonumber\\
&&+ \frac{2i}{3(1+\xi_{n}^{2})^{3}}
 \sum_{\alpha ,\beta <n}\Big(R^{\partial_{M}}_{i\alpha j\beta}(x_0)+R^{\partial_{M}}_{i\beta j\alpha}(x_0)\Big)
 \xi_{j}\xi_{\alpha}\xi_{\beta}\nonumber\\
&&+\frac{i\sum_{k}\xi_{k}} {4(1+\xi_{n}^{2})^{2}} \sum_{s\neq t <n}R^{\partial_{M}}_{k\gamma st}(x_0)\overline{c}(\widetilde{e_{s}})\overline{c}(\widetilde{e_{t}}).\nonumber\\
\end{eqnarray}
By the relation of the Clifford action and ${\rm tr}{(AB)}={\rm tr }{(BA)}$, we have
\begin{eqnarray}\sum_{s\neq t}\text{trace}[\overline{c}(\widetilde{e_{s}})\overline{c}(\widetilde{e_{t}})](x_{0})=0.
\end{eqnarray}
So we have
\begin{eqnarray}
&&\text{trace}\Big[\partial_{\xi_{n}}\partial_{\xi'}^{\alpha}\pi_{\xi_{n}}^{+}\sigma_{-2}(D_{T}^{-2})
\partial_{x'}^{\alpha}\sigma_{-3}(D_{T}^{-2})\Big](x_{0})\nonumber\\
&&=\frac{8\xi_{n}-24i} {3(\xi_{n}-i)^{3}(1+\xi_{n}^{2})^{3}}
\sum_{\alpha ,\beta <n}\Big(R^{\partial_{M}}_{i\alpha j\beta}(x_0)+R^{\partial_{M}}_{i\beta j\alpha}(x_0)\Big)
\xi_{i}\xi_{j}\xi_{\alpha}\xi_{\beta}\nonumber\\
&&+\frac{30i-10\xi_{n}} {3(\xi_{n}-i)^{3}(1+\xi_{n}^{2})^{2}}\sum_{i <n}R^{\partial_{M}}_{i \gamma i k }(x_0) \xi_{\gamma} \xi_{k}.
\end{eqnarray}
Since $R^{\partial_{M}}_{i\alpha j\beta}(x_0)=-R^{\partial_{M}}_{i\beta j\alpha}(x_0)$, and
by (3.59) and some calculations, we obtain
\begin{eqnarray}
\text{ Case \ (8) }
&=&\int_{|\xi'|=1}\int_{-\infty}^{+\infty}\frac{30i-10\xi_{n}} {3(\xi_{n}-i)^{3}(1+\xi_{n}^{2})^{2}}\sum_{i <n}R^{\partial_{M}}_{i \gamma i k }(x_0) \xi_{\gamma} \xi_{k} \nonumber\\
&=&\frac{1}{9}s_{\partial_{M}}
\frac{ 2\pi i }{4!}\bigg[\frac{15i-5\xi_{n}}
{(\xi_{n}+i)^{2}}\bigg]^{(4)}\bigg|_{\xi_{n}=i}\Omega_{5}\texttt{d}x' \nonumber\\
&=&\frac{5 }{16}s_{\partial_{M}}\pi\Omega_{5}\texttt{d}x'.
\end{eqnarray}

\textbf{Case (9)}: \ $r=-2, \ \ell=-3, \ k=1, \ j=0, \ |\alpha|=0$

From (2.15) and the Leibniz rule, we obtain
\begin{eqnarray}
\text{ Case \ (9)}&=&-\frac{1}{2}\int_{|\xi'|=1}\int_{-\infty}^{+\infty}\sum_{|\alpha|=1}
\text{trace}\Big[\partial_{\xi_{n}}\pi_{\xi_{n}}^{+}\sigma_{-2}(D_{T}^{-2})
\partial_{\xi_{n}}\partial_{x_{n}}\sigma_{-3}(D_{T}^{-2})\Big](x_{0})\texttt{d}\xi_{n}\sigma(\xi')\texttt{d}x' \nonumber\\
  &=&\frac{1}{2}\int_{|\xi'|=1}\int_{-\infty}^{+\infty}\sum_{|\alpha|=1}
\text{trace}\Big[\partial_{\xi_{n}}^{2}\pi_{\xi_{n}}^{+}\sigma_{-2}(D_{T}^{-2})
\partial_{x_{n}}\sigma_{-3}(D_{T}^{-2})\Big](x_{0})\texttt{d}\xi_{n}\sigma(\xi')\texttt{d}x'.
\end{eqnarray}
From (3.43), we have
\begin{equation}
\partial_{\xi_{n}}^{2}\pi_{\xi_{n}}^{+}\sigma_{-2}(D_{T}^{-2})(x_{0})\Big|_{|\xi'|=1}
=\frac{-i}{(\xi_{n}-i)^{3}}.
\end{equation}
By Lemma 3.1, Lemma 3.3 and some calculations

\begin{eqnarray}
\partial_{x_{n}}\sigma_{-3}(D_{T}^{-2})(x_{0})\Big|_{|\xi'|=1}&=&\partial_{x_{n}}\sigma_{-3}(D^{-2})(x_{0})\Big|_{|\xi'|=1}
+\partial_{x_{n}}(2i|\xi|^{-4}a^{n}\xi_{n})(x_{0})\Big|_{|\xi'|=1}\nonumber\\
&+&\partial_{x_{n}}(2i|\xi|^{-4}a^{k}\xi_{k})(x_{0})\Big|_{|\xi'|=1}.
\end{eqnarray}
Furtheremore we have
\begin{eqnarray}
\partial_{x_{n}}\sigma_{-3}(D_{T}^{-2})(x_{0})\Big|_{|\xi'|=1}&=&
\frac{2ih'(0) }{(1+\xi_{n}^{2})^{3}}\Big( -\frac{1}{2}h'(0)\sum_{k<n}\xi_{k}c(\widetilde{e_{k}})c(\widetilde{e_{n}})+3h'(0)\xi_{n} \Big)\nonumber\\
&-&\frac{i}{(1+\xi_{n}^{2})^{2}}\Big(\xi_{n}\big( 3h''(0)-\frac{9}{2}(h'(0))^{2}\big)-2\xi_{k}\big( \frac{3}{8}(h'(0))^{2}-\frac{1}{4}h''(0)\big)\sum_{t<n }c(\widetilde{e_{n}})c(\widetilde{e_{t}})\nonumber\\
&-&\frac{1}{4}\xi_{n}\big((h'(0))^{2}-h''(0)\big)\sum_{s\neq t<n }c(\widetilde{e_{s}})c(\widetilde{e_{t}})\Big)
+\frac{\big(h'(0)\big)^{2}}{4(1+\xi_{n}^{2})^{2}} \sum_{s\neq t,k <n}\xi_{k}\overline{c}(\widetilde{e_{s}})\overline{c}(\widetilde{e_{t}})\nonumber\\
&+&\frac{2i}{(1+\xi_{n}^{2})^{2}}\sum_{k,t<n}\xi_{k}\Big(\frac{3}{8}\big(h'(0)\big)^{2}-\frac{1}{4}h''(0)\Big)\overline{c}(\widetilde{e_{n}})\overline{c}(\widetilde{e_{t}})\nonumber\\
&+&\frac{i\xi_{n}} {4(1+\xi_{n}^{2})^{2}} \sum_{t <n}\Big(\big(h'(0)\big)^{2}-h''(0)\Big)\overline{c}(\widetilde{e_{s}})\overline{c}(\widetilde{e_{t}}).
\end{eqnarray}
By the relation of the Clifford action and ${\rm tr}{(AB)}={\rm tr }{(BA)}$, we have
\begin{equation}
\sum_{t<n}\text{trace}[\overline{c}(\widetilde{e_{s}})\overline{c}(\widetilde{e_{t}})](x_{0})=\sum_{t<n}\text{trace}[\overline{c}(\widetilde{e_{t}})\overline{c}(\widetilde{e_{t}})](x_{0})=48.
\end{equation}
By (3.62),(3.64),(3.65), we obtain
\begin{eqnarray}
\text{trace}\Big[\partial_{\xi_{n}}^{2}\pi_{\xi_{n}}^{+}\sigma_{-2}(D_{T}^{-2})
\partial_{x_{n}}\sigma_{-3}(D_{T}^{-2})\Big](x_{0})\Big|_{|\xi'|=1}
&=&\frac{12\xi_{n}\Big(\big(h'(0)\big)^{2}-h''(0)\Big)} {(\xi_{n}-i)^{3}(1+\xi_{n}^{2})^{2}}\nonumber\\
&+&\frac{\big(h'(0)\big)^{2}(84\xi_{n}+36\xi_{n}^{3})}
 {(\xi_{n}-i)^{3}(1+\xi_{n}^{2})^{3}}
+\frac{-24\xi_{n}h''(0)}
 {(\xi_{n}-i)^{3}(1+\xi_{n}^{2})^{2}}.\nonumber\\
\end{eqnarray}

Since
\begin{equation}
\partial_{x_{n}}c(\xi')(x_{0})=\sum_{j<n}\partial_{x_{n}}\xi_{j}c(dx_{j})=\sum_{j<n,1\leq l\leq n-1}\xi_{j}\partial_{x_{n}}(\sqrt{h}<dx^{j},\widetilde{e}^{l}>_{\partial M})c(\widetilde{e}^{l}),
\end{equation}
and
\begin{equation}
\int_{|\xi'|=1}\xi_1\cdots\xi_{2q+1}\sigma(\xi')=0.
\end{equation}
Therefore
\begin{eqnarray}
\text{ Case \ (9) }
&=&\frac{1}{2}\big(h'(0)\big)^{2}\int_{|\xi'|=1}\int_{-\infty}^{+\infty}
\frac{(84\xi_{n}+36\xi_{n}^{3})} {(\xi_{n}-i)^{3}(1+\xi_{n}^{2})^{3}}\texttt{d}\xi_{n}\sigma(\xi')\texttt{d}x'\nonumber\\
&&+\frac{1}{2}h''(0)\int_{|\xi'|=1}\int_{-\infty}^{+\infty}
\frac{-24\xi_{n}} {(\xi_{n}-i)^{3}(1+\xi_{n}^{2})^{2}}\texttt{d}\xi_{n}\sigma(\xi')\texttt{d}x'\nonumber\\
&&+\frac{1}{2}\int_{|\xi'|=1}\int_{-\infty}^{+\infty}\frac{12\xi_{n}\Big(\big(h'(0)\big)^{2}-h''(0)\Big)} {(\xi_{n}-i)^{3}(1+\xi_{n}^{2})^{2}}\texttt{d}\xi_{n}\sigma(\xi')\texttt{d}x'\nonumber\\
&=&\frac{1}{2}\big(h'(0)\big)^{2}  \frac{ 2\pi i }{5!}
\bigg[\frac{84\xi_{n}+36\xi_{n}^{3}}{(\xi_{n}+i)^{3}}\bigg]^{(5)}\bigg|_{\xi_{n}=i}\Omega_{5}\texttt{d}x'
+\frac{1}{2}h''(0) \frac{ 2\pi i }{4!}\bigg[\frac{-24\xi_{n}}
{(\xi_{n}+i)^{2}}\bigg]^{(4)}\bigg|_{\xi_{n}=i}\Omega_{5}\texttt{d}x'  \nonumber\\
&&+6\Big(\big(h'(0)\big)^{2}-h''(0)\Big)\frac{ 2\pi i }{4!}\bigg[\frac{\xi_{n}}
{(\xi_{n}+i)^{2}}\bigg]^{(4)}\bigg|_{\xi_{n}=i}\Omega_{5}\texttt{d}x'  \nonumber\\
&=&\Big(\frac{9}{16}h''(0)-\frac{45}{16}\big(h'(0)\big)^{2}\Big)\pi\Omega_{5}\texttt{d}x'.
\end{eqnarray}

\textbf{Case (10)}: \ $r=-3, \ \ell=-2, \ k=0, \ j=1, \ |\alpha|=0$

From (2.15), we have
\begin{equation}
\text{ Case \ (10)}=-\frac{1}{2}\int_{|\xi'|=1}\int_{-\infty}^{+\infty}
\text{trace}\Big[\partial_{x_{n}}\pi_{\xi_{n}}^{+}\sigma_{-3}(D_{T}^{-2})
\partial_{\xi_{n}}^{2}\sigma_{-2}(D_{T}^{-2})\Big](x_{0})\texttt{d}\xi_{n}\sigma(\xi')\texttt{d}x'.
\end{equation}
By the Leibniz rule, trace property and "++" and "-~-" vanishing
after the integration over $\xi_n$ in \cite{FGLS}, then
\begin{eqnarray}
&&\int^{+\infty}_{-\infty}{\rm trace}
\Big[\partial_{x_{n}}\pi_{\xi_{n}}^{+}\sigma_{-3}(D_{T}^{-2})
\partial_{\xi_{n}}^{2}\sigma_{-2}(D_{T}^{-2})\Big]\texttt{d}\xi_n \nonumber\\
&&= \int^{+\infty}_{-\infty}{\rm trace}
\Big[\partial_{x_{n}}\sigma_{-3}(D_{T}^{-2})\partial_{\xi_{n}}^{2}\sigma_{-2}(D_{T}^{-2})\Big]\texttt{d}\xi_n
-\int^{+\infty}_{-\infty}{\rm trace}
\Big[\partial_{x_{n}}\sigma_{-3}(D_{T}^{-2})\partial_{\xi_{n}}^{2}\pi_{\xi_{n}}^{+}\sigma_{-2}(D_{T}^{-2})\Big]\texttt{d}\xi_n.\nonumber\\
\end{eqnarray}
By case(9), we obtain
\begin{eqnarray}
&&\frac{1}{2}\int_{|\xi'|=1}\int_{-\infty}^{+\infty}
\text{trace}\Big[\partial_{\xi_{n}}^{2}\pi_{\xi_{n}}^{+}\sigma_{-2}(D_{T}^{-2})
\partial_{x_{n}}\sigma_{-3}(D_{T}^{-2})\Big](x_{0})\texttt{d}\xi_{n}\sigma(\xi')\texttt{d}x'\nonumber\\
&&=\Big(\frac{9}{16}h''(0)-\frac{45}{16}\big(h'(0)\big)^{2}\Big)\pi\Omega_{5}\texttt{d}x'.
\end{eqnarray}
By Lemma 3.1, Lemma 3.3 and a simple computation we obtain
\begin{equation}
\partial_{\xi_{n}}^{2}\sigma_{-2}(D_{T}^{-2})(x_{0})\Big|_{|\xi'|=1}
=\frac{6\xi_{n}^{2}-2}{(1+\xi_{n}^{2})^{3}}.
\end{equation}
By (3.73) and (3.64), we obtain
\begin{eqnarray}
\text{trace}\Big[\partial_{x_{n}}\sigma_{-3}(D_{T}^{-2})\partial_{\xi_{n}}^{2}\sigma_{-2}(D_{T}^{-2})\Big](x_{0})
&=&(h'(0))^{2}\frac{(84i\xi_{n}+36i\xi_{n}^{3})(-2+6\xi_{n}^{2})} {(1+\xi_{n}^{2})^{6}}
+h''(0)\frac{24i\xi_{n}(2-6\xi_{n}^{2})} {(1+\xi_{n}^{2})^{5}}.\nonumber\\
&+&\frac{12i(6\xi_{n}^{3}-2\xi_{n})\Big(\big(h'(0)\big)^{2}-h''(0)\Big)} {(1+\xi_{n}^{2})^{5}}.
\end{eqnarray}
By simple calculation, we obtain
\begin{equation}
-\frac{1}{2}\int_{|\xi'|=1}\int_{-\infty}^{+\infty}\int^{+\infty}_{-\infty}{\rm trace}
\Big[\partial_{x_{n}}\sigma_{-3}(D_{T}^{-2})\partial_{\xi_{n}}^{2}\sigma_{-2}(D_{T}^{-2})\Big]\texttt{d}\xi_n=\frac{45i}{32}\Big(\big(h'(0)\big)^{2}-h''(0)\Big)\pi\Omega_{5}\texttt{d}x'.
\end{equation}

Therefore
\begin{equation}
\text{ Case \ (10)}=\Big((\frac{9}{16}-\frac{45i}{32})h''(0)+(\frac{45i}{32}-\frac{45}{16})\big(h'(0)\big)^{2}\Big)\pi\Omega_{5}\texttt{d}x'.
\end{equation}

\textbf{Case (11)}: \ $r=-3, \ \ell=-2, \ k=0, \ j=0, \ |\alpha|=1$

From (2.15), we have
\begin{equation}
\text{ Case \ (11)}=-\int_{|\xi'|=1}\int_{-\infty}^{+\infty}\sum_{|\alpha|=1}
\text{trace}\Big[\partial_{\xi'}^{\alpha}\pi_{\xi_{n}}^{+}\sigma_{-3}(D_{T}^{-2})
\partial_{x'}^{\alpha}\partial_{\xi_{n}}\sigma_{-2}(D_{T}^{-2})\Big](x_{0})\texttt{d}\xi_{n}\sigma(\xi')\texttt{d}x' .
\end{equation}
By Lemma 3.1 and Lemma 3.3, for $i<n$, we have
\begin{equation}
\partial_{x_i}\sigma_{-2}(D_{T}^{-2})(x_0)=\partial_{x_i}\big(|\xi|^{-2}\big)(x_0)=0.
\end{equation}
So Case (11) vanishes.

\textbf{Case (12)}: \ $r=-3, \ \ell=-2, \ k=1, \ j=0, \ |\alpha|=0$

From (2.15) and the Leibniz rule, we have
\begin{eqnarray}
\text{ Case \ (12)}&=&-\frac{1}{2}\int_{|\xi'|=1}\int_{-\infty}^{+\infty}
\text{trace}\Big[\partial_{x_{n}}\pi_{\xi_{n}}^{+}\sigma_{-3}(D_{T}^{-2})
\partial_{\xi_{n}}\partial_{x_{n}}\sigma_{-2}(D_{T}^{-2})\Big](x_{0})\texttt{d}\xi_{n}\sigma(\xi')\texttt{d}x'\nonumber\\
&=&\frac{1}{2}\int_{|\xi'|=1}\int_{-\infty}^{+\infty}
\text{trace}\Big[\pi_{\xi_{n}}^{+}\sigma_{-3}(D_{T}^{-2})
\partial_{\xi_{n}}^{2}\partial_{x_{n}}\sigma_{-2}(D_{T}^{-2})\Big](x_{0})\texttt{d}\xi_{n}\sigma(\xi')\texttt{d}x'.
\end{eqnarray}
By the Leibniz rule, trace property and "++" and "-~-" vanishing
after the integration over $\xi_n$ in \cite{FGLS}, then
\begin{eqnarray}
&&\int^{+\infty}_{-\infty}{\rm trace}
\Big[\pi_{\xi_{n}}^{+}\sigma_{-3}(D_{T}^{-2})
\partial_{\xi_{n}}^{2}\partial_{x_{n}}\sigma_{-2}(D_{T}^{-2})\Big]\texttt{d}\xi_n \nonumber\\
&&= \int^{+\infty}_{-\infty}{\rm trace}
\Big[\sigma_{-3}(D_{T}^{-2})\partial_{\xi_{n}}^{2}\partial_{x_{n}}\sigma_{-2}(D_{T}^{-2})\Big]\texttt{d}\xi_n
-\int^{+\infty}_{-\infty}{\rm trace}
\Big[\sigma_{-3}(D_{T}^{-2})\partial_{\xi_{n}}^{2}\partial_{x_{n}}\pi_{\xi_{n}}^{+}\sigma_{-2}(D_{T}^{-2})\Big]\texttt{d}\xi_n.\nonumber\\
\end{eqnarray}

Similar to Case(7), we get the second term
\begin{equation}
-\frac{1}{2}\int_{|\xi'|=1}\int^{+\infty}_{-\infty}{\rm trace}
\Big[\sigma_{-3}(D_{T}^{-2})\partial_{\xi_{n}}^{2}\partial_{x_{n}}\pi_{\xi_{n}}^{+}\sigma_{-2}(D_{T}^{-2})\Big]\texttt{d}\xi_n=\frac{21}{8}\big(h'(0)\big)^{2}\pi\Omega_{5}\texttt{d}x'.
\end{equation}

From some direct computations, we obtain
\begin{equation}
\partial_{\xi_{n}}^{2}\partial_{x_{n}}\sigma_{-2}(D_{T}^{-2})(x_0)|_{|\xi'|=1}
=\frac{4-20\xi_n^{2}}{(1+\xi_n^{2})^{4}}h'(0).
\end{equation}

By the relation of the Clifford action and ${\rm tr}{(AB)}={\rm tr }{(BA)}$, we have
\begin{equation}
\text{trace}[\overline{c}(\widetilde{e_{n}})\overline{c}(\widetilde{e_{k}})](x_{0})=0.
\end{equation}

By (3.50) and (3.83), we obtain
\begin{eqnarray}
&&\text{trace}\Big[\sigma_{-3}(D_{T}^{-2})\partial_{\xi_{n}}^{2}\partial_{x_{n}}\sigma_{-2}(D_{T}^{-2})\Big](x_{0})
=-(h'(0))^{2}\frac{20i\xi_{n}-88i\xi_{n}^{3}-60i\xi_{n}^{5}} {(1+\xi_{n}^{2})^{7}}.\nonumber\\
\end{eqnarray}
From some direct computations, we obtain
\begin{equation}
\int_{-\infty}^{+\infty}\frac{-20i\xi_{n}+88i\xi_{n}^{3}+60i\xi_{n}^{5}} {(1+\xi_{n}^{2})^{7}}\texttt{d}\xi_{n}
=\frac{ 2\pi i }{6!}\bigg[\frac{-20i\xi_{n}+88i\xi_{n}^{3}+60i\xi_{n}^{5}}{(\xi_{n}+i)^{7}}\bigg]^{(6)}\Big|_{\xi_{n}=i}=0.
\end{equation}
Therefore the first term of (3.80) is vanishes. So
\begin{equation}
\text{ Case \ (12)}=\frac{21}{8}\big(h'(0)\big)^{2}\pi\Omega_{5}\texttt{d}x'.
\end{equation}
\textbf{Case (13)}: \ $r=-3, \ \ell=-3, \ k=0, \ j=0, \ |\alpha|=0$

From (2.15) and the Leibniz rule, we have
\begin{eqnarray}
\text{ Case \ (13)}&=&-i\int_{|\xi'|=1}\int_{-\infty}^{+\infty}
\text{trace}\Big[\pi_{\xi_{n}}^{+}\sigma_{-3}(D_{T}^{-2})
\partial_{\xi_{n}}\sigma_{-3}(D_{T}^{-2})\Big](x_{0})\texttt{d}\xi_{n}\sigma(\xi')\texttt{d}x'\nonumber\\
&=&i\int_{|\xi'|=1}\int_{-\infty}^{+\infty}
\text{trace}\Big[\partial_{\xi_{n}}\pi_{\xi_{n}}^{+}\sigma_{-3}(D_{T}^{-2})
\sigma_{-3}(D_{T}^{-2})\Big](x_{0})\texttt{d}\xi_{n}\sigma(\xi')\texttt{d}x'.
\end{eqnarray}

By (2.13), we obtain
\begin{equation}
\pi^+_{\xi_n}\left[\frac{-5i\xi_n-3i\xi_n^3}{(1+\xi_n^2)^3}\right]
=\frac{9i-7\xi_n}{8(\xi_n-i)^{3}}.
\end{equation}

Then we obtain
\begin{eqnarray}
\partial_{\xi_{n}}\pi_{\xi_{n}}^{+}\sigma_{-3}(D_{T}^{-2})(x_0)|_{|\xi'|=1}&=&\frac{3i-\xi_n}{4(\xi_n-i)^3}h'(0)\sum_{k<n}\xi_k
c(\widetilde{e_k})c(\widetilde{e_n})+h'(0)\frac{7\xi_n-10i}{4(\xi_n-i)^4}\nonumber\\
&-&\frac{\xi_n^{2}-4\xi_n+3i}{8(\xi_n-i)^4}ih'(0)\sum_{k<n}\xi_k
\overline{c}(\widetilde{e_k})\overline{c}(\widetilde{e_n}).
\end{eqnarray}

By the relation of the Clifford action and ${\rm tr}{(AB)}={\rm tr }{(BA)}$, we have
\begin{eqnarray}
&&\sum_{k,l<n}\xi_l\xi_k\text{trace}[\overline{c}(\widetilde{e_{n}})\overline{c}(\widetilde{e_{k}})\overline{c}(\widetilde{e_{n}})\overline{c}(\widetilde{e_{l}})](x_{0})=-8\sum_{k<n}\xi_k^{2},\nonumber\\
&&\sum_{k,l<n}\text{trace}[\overline{c}(\widetilde{e_{n}})\overline{c}(\widetilde{e_{k}})c(\widetilde{e_{n}})c(\widetilde{e_{l}})](x_{0})=0.\nonumber\\
\end{eqnarray}

Then by (3.89),(3.90) we have
\begin{eqnarray}
\text{trace}\Big[\partial_{\xi_{n}}\pi_{\xi_{n}}^{+}\sigma_{-3}(D_{T}^{-2})\sigma_{-3}(D_{T}^{-2})\Big](x_{0})
&=&(h'(0))^{2}\frac{2(10i-7\xi_n)(5i\xi_n+3i\xi_n^{3})}{(\xi_n-i)^4(1+\xi_n^2)^3}\nonumber\\
&-&(h'(0))^{2}\frac{(i\xi_n^{2}+4\xi_n-3i)}{2(\xi_n-i)^4(1+\xi_n^2)^2}\sum_{k<n}\xi_k^{2}.\nonumber\\
\end{eqnarray}
Therefore
\begin{eqnarray}
\text{ Case \ (13) }
&=&i\big(h'(0)\big)^{2}\int_{|\xi'|=1}\int_{-\infty}^{+\infty}
\text{trace}\Big[\partial_{\xi_{n}}\pi_{\xi_{n}}^{+}\sigma_{-3}(D_{T}^{-2})\sigma_{-3}(D_{T}^{-2})\Big](x_{0})\texttt{d}\xi_{n}\sigma(\xi')\texttt{d}x'\nonumber\\
&=&i\big(h'(0)\big)^{2}  \frac{ 2\pi i }{6!}
\bigg[\frac{-3i-96\xi_n-72i\xi_n^{2}-56\xi_n^{3}-41i\xi_n^{4}}
{(\xi_{n}+i)^{3}}\bigg]^{(6)}\bigg|_{\xi_{n}=i}\Omega_{5}\texttt{d}x'  \nonumber\\
&+&\frac{1}{2}\big(h'(0)\big)^{2}  \frac{ 2\pi i }{5!}
\bigg[\frac{\xi_n^{2}-4i\xi_n-3)}
{(\xi_{n}+i)^{2}}\bigg]^{(5)}\bigg|_{\xi_{n}=i}\Omega_{5}\texttt{d}x'  \nonumber\\
&=&(-\frac{5}{32}-\frac{57}{8} )(h'(0))^{2}\pi\Omega_{5}dx'.
\end{eqnarray}

\textbf{Case (14)}: \ $r=-2, \ \ell=-4, \ k=0, \ j=0, \ |\alpha|=0$

From (2.15) and the Leibniz rule, we have
\begin{eqnarray}
\text{ Case \ (14)}&=&-i\int_{|\xi'|=1}\int_{-\infty}^{+\infty}
\text{trace}\Big[\pi_{\xi_{n}}^{+}\sigma_{-2}(D_{T}^{-2})
\partial_{\xi_{n}}\sigma_{-4}(D_{T}^{-2})\Big](x_{0})\texttt{d}\xi_{n}\sigma(\xi')\texttt{d}x'\nonumber\\
&=&i\int_{|\xi'|=1}\int_{-\infty}^{+\infty}
\text{trace}\Big[\partial_{\xi_{n}}\pi_{\xi_{n}}^{+}\sigma_{-2}(D_{T}^{-2})
\sigma_{-4}(D_{T}^{-2})\Big](x_{0})\texttt{d}\xi_{n}\sigma(\xi')\texttt{d}x'.
\end{eqnarray}
From (3.43), we have
\begin{equation}
\partial_{\xi_{n}}\pi_{\xi_{n}}^{+}\sigma_{-2}(D_{T}^{-2})(x_{0})\Big|_{|\xi'|=1}
=\frac{i}{2(\xi_{n}-i)^{2}}.
\end{equation}
By Lemma 3.2 and some direct calculation we have
\begin{eqnarray}
\sigma_{-4}(D_{T}^{-2})(x_{0})\Big|_{|\xi'|=1}&=&\sigma_{-4}(D^{-2})(x_{0})\Big|_{|\xi'|=1}-12\frac{h'(0)}{(1+\xi_{n}^{2})^{3}}+\frac{2\big(h'(0)\big)^{2}}{(1+\xi_{n}^{2})^{4}}\sum_{k<n,l}\xi_{k}\xi_{l}c(\widetilde{e_n})c(\widetilde{e_k})\nonumber\\
&-&\frac{7\big(h'(0)\big)^{2}}{16(1+\xi_{n}^{2})^{2}}\sum_{k,l<n}\xi_{k}\xi_{l}\overline{c}(\widetilde{e_n})\overline{c}(\widetilde{e_k})\overline{c}(\widetilde{e_n})\overline{c}(\widetilde{e_l})
+\frac{21\big(h'(0)\big)^{2}}{4(1+\xi_{n}^{2})^{2}}\sum_{k<n}\overline{c}(\widetilde{e_n})\overline{c}(\widetilde{e_k})\nonumber\\
&+&\frac{1}{8(1+\xi_{n}^{2})^{2}}\sum_{ijkl}R^{\partial_{M}}_{ijkl}(x_0)\overline{c}(\widetilde{e_i})\overline{c}(\widetilde{e_j})c(\widetilde{e_k})c(\widetilde{e_l})
-\frac{T\sum_{i}c(\widetilde{e_i})\bar{c}(\nabla_{\widetilde{e_i}}^{TM}V)+T^{2}|V|^{2}}{(1+\xi_{n}^{2})^{2}}\nonumber\\
&+&\frac{2h'(0)\sum_{k}\xi_{k}}{(1+\xi_{n}^{2})^{3}}-\frac{(\sum_{k<n}\xi_{n}\xi_{k})+\xi_{n}^{2}}{2(1+\xi_{n}^{2})}\sum_{t<n,s}\big((h'(0))^{2}-h''(0)\big)\overline{c}(\widetilde{e_s})\overline{c}(\widetilde{e_t})\nonumber\\
&-&4\frac{\sum_{k,l<n}\xi_{k}\xi_{l}+\sum_{k<n}\xi_{k}\xi_{n}}{(1+\xi_{n}^{2})}\sum_{t<n,s}\Big(\frac{3}{8}\big(h'(0)\big)^{2}-\frac{1}{4}h''(0)\Big)\overline{c}(\widetilde{e_s})\overline{c}(\widetilde{e_t})\nonumber\\
&-&\frac{\big(h'(0)\big)^{2}}{4(1+\xi_{n}^{2})^{3}}\sum_{k,l<n}\xi_{k}\xi_{l}\overline{c}(\widetilde{e_n})\overline{c}(\widetilde{e_k})\overline{c}(\widetilde{e_n})\overline{c}(\widetilde{e_l})-\frac{7\big(h'(0)\big)^{2}}{8(1+\xi_{n}^{2})^{2}}\overline{c}(\widetilde{e_n})\overline{c}(\widetilde{e_k})c(\widetilde{e_n})c(\widetilde{e_l})\nonumber\\
&-&\frac{\big(h'(0)\big)^{2}}{8(1+\xi_{n}^{2})^{3}}\sum_{k,l<n}\xi_{k}\xi_{l}c(\widetilde{e_n})c(\widetilde{e_k})\overline{c}(\widetilde{e_n})\overline{c}(\widetilde{e_l}).\nonumber\\
\end{eqnarray}
Where
\begin{eqnarray}
\sigma_{-4}(D^{-2})(x_{0})\Big|_{|\xi'|=1}
&=&\frac{-\big(h'(0)\big)^{2}}{4(1+\xi_{n}^{2})^{3}}c(\widetilde{e_k})c(\widetilde{e_n})c(\widetilde{e_l})c(\widetilde{e_n})
- \frac{9\big(h'(0)\big)^{2}}{(1+\xi_{n}^{2})^{3}}\xi_{n}^{3}\xi_{\mu}\xi_{l}\nonumber\\
&&+\frac{\big(h'(0)\big)^{2}}{4(1+\xi_{n}^{2})^{2}}\xi_{k}\xi_{l}c(\widetilde{e_k})c(\widetilde{e_n})c(\widetilde{e_l})c(\widetilde{e_n})
-\frac{1}{4(1+\xi_{n}^{2})^{2}}s(x_{0})\nonumber\\
&&-\frac{5}{3(1+\xi_{n}^{2})^{3}}\xi_{k}\xi_{l}\sum_{i<n}R^{\partial_{M}}_{ik il}(x_{0})
-\frac{6}{(1+\xi_{n}^{2})^{3}}h''(0)\xi_{n}^{2}\nonumber\\
&&-\frac{4}{3(1+\xi_{n}^{2})^{4}}\xi_{k}\xi_{l}\xi_{\gamma}\xi_{\delta}
 \sum_{\gamma ,\delta <n}\Big(R^{\partial_{M}}_{k\gamma l\delta}(x_0)+R^{\partial_{M}}_{l\gamma k\delta}(x_0)\Big)
 +\frac{4h''(0)}{(1+\xi_{n}^{2})^{4}}\xi_{n}^{2}\nonumber\\
 &&-\frac{1}{3(1+\xi_{n}^{2})^{3}}\xi_{\alpha}\xi_{\beta}
 \sum_{\alpha ,\beta <n}\Big(R^{\partial_{M}}_{k\alpha l\beta}(x_0)+R^{\partial_{M}}_{l\beta k\alpha}(x_0)\Big)
 +\frac{h''(0)}{(1+\xi_{n}^{2})^{3}}\nonumber\\
 &&+\frac{2+3\xi_{n}+10\xi_{n}^{2}+12\xi_{n}^{3}-4\xi_{n}^{4}+9\xi_{n}^{5}}
{(1+\xi_{n}^{2})^{5}}\big(h'(0)\big)^{2}.
\end{eqnarray}

By the relation of the Clifford action and ${\rm tr}{(AB)}={\rm tr }{(BA)}$, we have
\begin{eqnarray}
&&{\rm trace}\big[c(\widetilde{e_k})c(\widetilde{e_n})c(\widetilde{e_l})c(\widetilde{e_n})\big],\nonumber\\
&&{\rm trace}\big[c(\widetilde{e_k})c(\widetilde{e_n})c(\widetilde{e_l})c(\widetilde{e_n})\big],\nonumber\\
&&{\rm trace}\big[c(\widetilde{e_i})\bar{c}(\nabla_{\widetilde{e_i}}^{TM}V)\big]=0.\nonumber\\
\end{eqnarray}
By (3.94)-(3.97), we obtain
\begin{eqnarray}
&&\text{trace}\Big[\partial_{\xi_{n}}\pi_{\xi_{n}}^{+}\sigma_{-2}(D_{T}^{-2})
\sigma_{-4}(D_{T}^{-2})\Big](x_{0})=\text{trace}\Big[\partial_{\xi_{n}}\pi_{\xi_{n}}^{+}\sigma_{-2}(D^{-2})
\sigma_{-4}(D^{-2})\Big](x_{0})-\frac{48ih'(0)}{(\xi_{n}-i)^{2}(1+\xi_{n}^{2})^{3}}\nonumber\\
&&-\frac{2i(h'(0))^{2}}{(\xi_{n}-i)^{2}(1+\xi_{n}^{2})^{3}}+\frac{21i\big((h'(0))^{2}-h''(0)\big)}{(\xi_{n}-i)^{2}(1+\xi_{n}^{2})^{2}}+\frac{is_{M}}{2(\xi_{n}-i)^{2}(1+\xi_{n}^{2})^{2}}\nonumber\\
&&+\frac{21i(h'(0))^{2}}{2(\xi_{n}-i)^{2}(1+\xi_{n}^{2})^{2}}+\frac{-i4T^{2}|V|^{2}}{(\xi_{n}-i)^{2}(1+\xi_{n}^{2})^{2}}+\frac{8ih'(0)\xi_{n}}{(\xi_{n}-i)^{2}(1+\xi_{n}^{2})^{3}}\nonumber\\
&&\nonumber\\
&&-\frac{\big(12i\sum_{k<n}\xi_{k}\xi_{n}+12i\xi_{n}^{2}\big)\big((h'(0))^{2}-h''(0)\big)}{(\xi_{n}-i)^{2}(1+\xi_{n}^{2})}\nonumber\\
&&-\frac{\big(96i\sum_{k<n}\xi_{k}\xi_{n}+96i\sum_{k<n}\xi_{k}^{2}\big)\big(\frac{3}{8}\big(h'(0)\big)^{2}-\frac{1}{4}h''(0)\big)}{(\xi_{n}-i)^{2}(1+\xi_{n}^{2})}.\nonumber\\
\end{eqnarray}
where
\begin{eqnarray}
&&\text{trace}\Big[\partial_{\xi_{n}}\pi_{\xi_{n}}^{+}\sigma_{-2}(D^{-2})
\sigma_{-4}(D^{-2})\Big](x_{0})\Big|_{|\xi'|=1}\nonumber\\
&&=\frac{-i\big(h'(0)\big)^{2}}{8(\xi_{n}-i)^{2}(1+\xi_{n}^{2})^{3}}
{\rm tr}\big[c(\widetilde{e_k})c(\widetilde{e_n})c(\widetilde{e_l})c(\widetilde{e_n})\big]
- \frac{36i\big(h'(0)\big)^{2}\xi_{n}^{3}\sum_{k,l<n}\xi_{k}\xi_{l}}{(\xi_{n}-i)^{2}(1+\xi_{n}^{2})^{3}}\nonumber\\
&&+\frac{i\big(h'(0)\big)^{2}\sum_{k,l<n}\xi_{k}\xi_{l}}{8(\xi_{n}-i)^{2}(1+\xi_{n}^{2})^{2}}
{\rm tr}\big[c(\widetilde{e_k})c(\widetilde{e_n})c(\widetilde{e_l})c(\widetilde{e_n})\big]
-\frac{is(x_{0})}{(\xi_{n}-i)^{2}(1+\xi_{n}^{2})^{2}}\nonumber\\
&&-\frac{20i\sum_{k,l<n}\xi_{k}\xi_{l}\sum_{i<n}R^{\partial_{M}}_{ik il}(x_{0})}{3(\xi_{n}-i)^{2}(1+\xi_{n}^{2})^{3}}
-\frac{24ih''(0)\xi_{n}^{2}}{(\xi_{n}-i)^{2}(1+\xi_{n}^{2})^{3}}\nonumber\\
&&-\frac{16i\xi_{k}\xi_{l}\xi_{\gamma}\xi_{\delta}
 \sum_{\gamma ,\delta <n}\Big(R^{\partial_{M}}_{k\gamma l\delta}(x_0)+R^{\partial_{M}}_{l\gamma k\delta}(x_0)\Big)}{3(\xi_{n}-i)^{2}(1+\xi_{n}^{2})^{4}}
 \nonumber\\
 &&-\frac{4i\sum_{k,l<n}\xi_{k}\xi_{l}\sum_{\alpha ,\beta <n}\Big(R^{\partial_{M}}_{k\alpha l\beta}(x_0)+R^{\partial_{M}}_{l\beta k\alpha}(x_0)\Big)}{3(\xi_{n}-i)^{2}(1+\xi_{n}^{2})^{3}}
 \nonumber\\
 &&+\frac{16ih''(0)\xi_{n}^{2}}{(\xi_{n}-i)^{2}(1+\xi_{n}^{2})^{4}}
 +\frac{4ih''(0)}{(\xi_{n}-i)^{2}(1+\xi_{n}^{2})^{3}}\nonumber\\
  &&+\frac{4i\big(h'(0)\big)^{2}(2+3\xi_{n}+10\xi_{n}^{2}+12\xi_{n}^{3}-4\xi_{n}^{4}+9\xi_{n}^{5})}
{(\xi_{n}-i)^{2}(1+\xi_{n}^{2})^{5}}.
\end{eqnarray}
Therefore
\begin{eqnarray}
\text{ Case \ (14) }&=&i\int_{|\xi'|=1}\int_{-\infty}^{+\infty}
\text{trace}\Big[\partial_{\xi_{n}}\pi_{\xi_{n}}^{+}\sigma_{-2}(D_{T}^{-2})
\sigma_{-4}(D_{T}^{-2})\Big](x_{0})\texttt{d}\xi_{n}\sigma(\xi')\texttt{d}x'\nonumber\\
&=&s_{M}(x_{0})\frac{2\pi i}{3!}\bigg[\frac{1}{(\xi_{n}+i)^{2}}\bigg]^{(3)}\bigg|_{\xi_{n}=i}\Omega_{5}\texttt{d}x'+s_{\partial_{M}}(x_{0})\frac{2\pi i}{4!}\bigg[\frac{14}
{9(\xi_{n}+i)^{3}}\bigg]^{(4)}\bigg|_{\xi_{n}=i}\Omega_{5}\texttt{d}x'
\nonumber\\
&&+\big(h'(0)\big)^{2}  \frac{ 2\pi i }{6!}
\bigg[\frac{-53-72\xi_{n}-33\xi_{n}^{2}-288\xi_{n}^{3}+525\xi_{n}^{4}-216\xi_{n}^{5}+217\xi_{n}^{6}}
{6(\xi_{n}+i)^{5}}\bigg]^{(6)}\bigg|_{\xi_{n}=i}\Omega_{5}\texttt{d}x'  \nonumber\\
&&+ h''(0) \frac{ 2\pi i }{5!}\bigg[\frac{-4(1-\xi_{n}^{2}-6\xi_{n}^{4})}
{(\xi_{n}+i)^{4}}\bigg]^{(5)}\bigg|_{\xi_{n}=i}\Omega_{5}\texttt{d}x' \nonumber\\
&&+\big( 48h'(0)+2(h'(0))^{2}\big)\frac{ 2\pi i }{4!}\bigg[\frac{1}
{(\xi_{n}+i)^{3}}\bigg]^{(4)}\bigg|_{\xi_{n}=i}\Omega_{5}\texttt{d}x'\nonumber\\
&&+\big(4T^{2}|V|^{2}-21\big((h'(0))^{2}-h''(0)\big)-\frac{1}{2}s_{M}-\frac{21}{2}(h'(0))^{2}\big)\frac{ 2\pi i }{3!}\bigg[\frac{1}
{(\xi_{n}+i)^{2}}\bigg]^{(3)}\bigg|_{\xi_{n}=i}\Omega_{5}\texttt{d}x' \nonumber\\
&&-8h'(0)\frac{ 2\pi i }{4!}\bigg[\frac{\xi_{n}}
{(\xi_{n}+i)^{3}}\bigg]^{(4)}\bigg|_{\xi_{n}=i}\Omega_{5}\texttt{d}x'+12\big((h'(0))^{2}-h''(0)\big)\frac{ 2\pi i }{2!}\bigg[\frac{\xi_{n}^{2}}
{(\xi_{n}+i)}\bigg]^{(2)}\bigg|_{\xi_{n}=i}\Omega_{5}\texttt{d}x'\nonumber\\
&&+96\big(\frac{3}{8}\big(h'(0)\big)^{2}-\frac{1}{4}h''(0)\big)\frac{ 2\pi i }{2!}\bigg[\frac{1}
{(\xi_{n}+i)}\bigg]^{(2)}\bigg|_{\xi_{n}=i}\Omega_{5}\texttt{d}x'\nonumber\\
&=&\Big(\frac{-1}{4}s_{M}(x_{0})-(\frac{45}{4}+\frac{5i}{8})h'(0)-\big(\frac{23}{12}+\frac{3i}{2}\big)\big(h'(0)\big)^{2}
+\frac{235}{64}h''(0)+\frac{47}{96}s_{\partial_{M}}(x_{0})-T^{2}|V|^{2}\Big)\pi\Omega_{5}\texttt{d}x'.\nonumber\\
\end{eqnarray}

\textbf{Case (15)}: \ $r=-4, \ \ell=-2, \ k=0, \ j=0, \ |\alpha|=0$

From (2.15), we have
\begin{equation}
\text{ Case \ (15)}=-i\int_{|\xi'|=1}\int_{-\infty}^{+\infty}
\text{trace}\Big[\pi_{\xi_{n}}^{+}\sigma_{-4}(D_{T}^{-2})
\partial_{\xi_{n}}\sigma_{-2}(D_{T}^{-2})\Big](x_{0})\texttt{d}\xi_{n}\sigma(\xi')\texttt{d}x'.
\end{equation}
By the Leibniz rule, trace property and "++" and "-~-" vanishing after the integration over $\xi_n$ in \cite{FGLS}, then
\begin{eqnarray}
&&\int^{+\infty}_{-\infty}{\rm trace}
\Big[\pi_{\xi_{n}}^{+}\sigma_{-4}(D_{T}^{-2})
\partial_{\xi_{n}}\sigma_{-2}(D_{T}^{-2})\Big]\texttt{d}\xi_n \nonumber\\
&&= \int^{+\infty}_{-\infty}{\rm trace}
\Big[\sigma_{-4}(D_{T}^{-2})
\partial_{\xi_{n}}\sigma_{-2}(D_{T}^{-2})\Big]\texttt{d}\xi_n
-\int^{+\infty}_{-\infty}{\rm trace}
\Big[\sigma_{-4}(D_{T}^{-2})
\partial_{\xi_{n}}\pi_{\xi_{n}}^{+}\sigma_{-2}(D_{T}^{-2})\Big]\texttt{d}\xi_n.\nonumber\\
\end{eqnarray}
By Case (14), we obtain
\begin{eqnarray}
&&i\int_{|\xi'|=1}\int_{-\infty}^{+\infty}{\rm trace}
\Big[\sigma_{-4}(D_{T}^{-2})
\partial_{\xi_{n}}\pi_{\xi_{n}}^{+}\sigma_{-2}(D_{T}^{-2})\Big]\texttt{d}\xi_n\sigma(\xi')\texttt{d}x'\nonumber\\
&&=\Big(\frac{-1}{4}s_{M}(x_{0})-(\frac{45}{4}+\frac{5i}{8})h'(0)-\big(\frac{23}{12}+\frac{3i}{2}\big)\big(h'(0)\big)^{2}
+\frac{235}{64}h''(0)+\frac{47}{96}s_{\partial_{M}}(x_{0})-T^{2}|V|^{2}\Big)\pi\Omega_{5}\texttt{d}x'\nonumber\\
\end{eqnarray}
By Lemma 3.1 and Lemma 3.3, we obtain
\begin{equation}
\partial_{\xi_{n}}\sigma_{-2}(D_{T}^{-2})(x_{0})\Big|_{|\xi'|=1}
=\frac{-2\xi_{n}}{(1+\xi_{n}^{2})^{2}}.
\end{equation}
By (3.95),(3.96), and (3.104), we obtain
\begin{eqnarray}
&&\text{trace}\Big[\partial_{\xi_{n}}\pi_{\xi_{n}}^{+}\sigma_{-2}(D_{T}^{-2})
\sigma_{-4}(D_{T}^{-2})\Big](x_{0})=\text{trace}\Big[\partial_{\xi_{n}}\pi_{\xi_{n}}^{+}\sigma_{-2}(D^{-2})
\sigma_{-4}(D^{-2})\Big](x_{0})+\frac{192h'(0)\xi_{n}}{(1+\xi_{n}^{2})^{5}}\nonumber\\
&&-\frac{16(h'(0))^{2}\xi_{n}}{(1+\xi_{n}^{2})^{5}}+\frac{12\big((h'(0))^{2}-h''(0)\big)\xi_{n}}{(1+\xi_{n}^{2})^{4}}-\frac{2s_{M}\xi_{n}}{(1+\xi_{n}^{2})^{4}}\nonumber\\
&&+\frac{12(h'(0))^{2}\xi_{n}}{(1+\xi_{n}^{2})^{2}}-\frac{16T^{2}|V|^{2}}{(1+\xi_{n}^{2})^{4}}-\frac{32h'(0)\sum_{k}\xi_{k}\xi_{n}}{(1+\xi_{n}^{2})^{3}}\nonumber\\
&&-\frac{48\sum_{k<n}\xi_{k}\xi^{2}_{n}\big((h'(0))^{2}-h''(0)\big)}{(1+\xi_{n}^{2})^{5}}-\frac{48\xi_{n}^{2}\big((h'(0))^{2}-h''(0)\big)\xi_{n}^{3}}{(1+\xi_{n}^{2})^{4}}\nonumber\\
&&-\frac{-384\sum_{k<n}\xi_{k}\xi_{n}^{2}\Big(\frac{3}{8}\big(h'(0)\big)^{2}-\frac{1}{4}h''(0)\Big)}{(1+\xi_{n}^{2})^{3}}-\frac{384\sum_{k<n}\xi_{k}^{2}\Big(\frac{3}{8}\big(h'(0)\big)^{2}-\frac{1}{4}h''(0)\Big)\xi_{n}}{(1+\xi_{n}^{2})^{3}},\nonumber\\
\end{eqnarray}
where
\begin{eqnarray}
&&\text{trace}\Big[\sigma_{-4}(D^{-2})
\partial_{\xi_{n}}\sigma_{-2}(D^{-2})\Big](x_{0})\Big|_{|\xi'|=1}\nonumber\\
&&=\frac{\xi_{n}\big(h'(0)\big)^{2}}{2(1+\xi_{n}^{2})^{5}}
{\rm tr}\big[c(\widetilde{e_k})c(\widetilde{e_n})c(\widetilde{e_l})c(\widetilde{e_n})\big]
+\frac{144\big(h'(0)\big)^{2}\xi_{n}^{4}\sum_{k,l<n}\xi_{k}\xi_{l}}{(1+\xi_{n}^{2})^{5}}\nonumber\\
&&-\frac{\big(h'(0)\big)^{2}\sum_{k,l<n}\xi_{k}\xi_{l}}{2(1+\xi_{n}^{2})^{4}}
{\rm tr}\big[c(\widetilde{e_k})c(\widetilde{e_n})c(\widetilde{e_l})c(\widetilde{e_n})\big]
-\frac{4s(x_{0})\xi_{n}}{(1+\xi_{n}^{2})^{4}}\nonumber\\
&&+\frac{80\xi_{n}\sum_{k,l<n}\xi_{k}\xi_{l}\sum_{i<n}R^{\partial_{M}}_{ikil}(x_{0})}{3(1+\xi_{n}^{2})^{5}}
+\frac{96h''(0)\xi_{n}^{3}}{(1+\xi_{n}^{2})^{5}}\nonumber\\
&&-\frac{64\xi_{n}\sum_{k,l,\gamma,\delta<n}\xi_{k}\xi_{l}\xi_{\gamma}\xi_{\delta}}{3(1+\xi_{n}^{2})^{6}}
 \sum_{\gamma ,\delta <n}\Big(R^{\partial_{M}}_{k\gamma l\delta}(x_0)+R^{\partial_{M}}_{l\gamma k\delta}(x_0)\Big)
 \nonumber\\
 &&+\frac{16\xi_{n}\sum_{k,l<n}\xi_{k}\xi_{l}}{3(1+\xi_{n}^{2})^{5}}
 \sum_{\alpha ,\beta <n}\Big(R^{\partial_{M}}_{k\alpha l\beta}(x_0)+R^{\partial_{M}}_{l\beta k\alpha}(x_0)\Big)
 \nonumber\\
 &&-\frac{64\xi_{n}^{3}h''(0)}{(1+\xi_{n}^{2})^{6}}
 -\frac{18\xi_{n}h''(0)}{(1+\xi_{n}^{2})^{5}}\nonumber\\
  &&-\frac{16\xi_{n}\big(h'(0)\big)^{2}(2+3\xi_{n}+10\xi_{n}^{2}+12\xi_{n}^{3}-4\xi_{n}^{4}+9\xi_{n}^{5})}
{(1+\xi_{n}^{2})^{7}}.
\end{eqnarray}

By similar calculations, we get
\begin{equation}
-i\int_{|\xi'|=1}\int^{+\infty}_{-\infty}{\rm trace}
\Big(\sigma_{-4}(D_{T}^{-2})
\partial_{\xi_{n}}\sigma_{-2}(D_{T}^{-2})\Big)\texttt{d}\xi_n=3i(h'(0))^{2}\pi\Omega_{5}\texttt{d}x'.
\end{equation}
Therefore
\begin{eqnarray}
&&\text{ Case \ (15)}=
\Big(\frac{-1}{4}s_{M}(x_{0})-(\frac{45}{4}+\frac{5i}{8})h'(0)-\big(\frac{23}{12}-\frac{3i}{2}\big)\big(h'(0)\big)^{2}
+\frac{235}{64}h''(0)+\frac{47}{96}s_{\partial_{M}}(x_{0})-T^{2}|V|^{2}\Big)\pi\Omega_{5}\texttt{d}x'.\nonumber\\
\end{eqnarray}

Now $\Phi$  is the sum of the \textbf{case ($1,2,\cdots,15$)}, so
\begin{eqnarray}
\Phi=\sum_{I=1}^{15} \textbf{case I}&=&\Big(\frac{-1}{2}s_{M}(x_{0})+\frac{35}{24}s_{\partial_{M}}(x_{0})-2T^{2}|V|^{2}-(\frac{45}{2}+\frac{5i}{4})h'(0)\nonumber\\
&&+\big(\frac{45i}{32}-\frac{3947}{384}\big)\big(h'(0)\big)^{2}
+\big(\frac{247}{32}-\frac{45i}{32}\big)h''(0)\Big)\pi\Omega_{5}\texttt{d}x'.
\end{eqnarray}
Hence we have the conclusion as follows.

\begin{thm}
 Let M be a $7$-dimensional spin compact manifold with the boundary $\partial M$. Then we get the volumes associated to Witten deformation $D_{T}$ on $\widehat{M}$
 \begin{eqnarray}
 \widetilde{{\rm Wres}}[\pi^+D_{T}^{-2}\circ\pi^+D_{T}^{-2}]&=&\int_{\partial M}\Big(\frac{-1}{2}s_{M}(x_{0})+\frac{35}{24}s_{\partial_{M}}(x_{0})-2T^{2}|V|^{2}-(\frac{45}{2}+\frac{5i}{4})h'(0)\nonumber\\
&&+\big(\frac{45i}{32}-\frac{3947}{384}\big)\big(h'(0)\big)^{2}
+\big(\frac{247}{32}-\frac{45i}{32}\big)h''(0)\Big)\pi\Omega_{5}\texttt{d}x'.
\end{eqnarray}
\end{thm}
\section{The gravitational action for 7-dimensional manifolds with boundary}
 Firstly, we recall the Einstein-Hilbert action for manifolds with boundary (see \cite{Wa3} or \cite{Wa4}),
 \begin{equation}
I_{\rm Gr}=\frac{1}{16\pi}\int_Ms{\rm dvol}_M+2\int_{\partial M}K{\rm dvol}_{\partial_M}:=I_{\rm {Gr,i}}+I_{\rm {Gr,b}},
\end{equation}
 where
 \begin{equation}
 K=\sum_{1\leq i,j\leq {n-1}}K_{i,j}g_{\partial M}^{i,j};~~K_{i,j}=-\Gamma^n_{i,j},
\end{equation}
 and $K_{i,j}$ is the second fundamental form, or extrinsic curvature. Taking the metric in Section 2, then for $n=7$ we have
 \begin{equation}
K(x_0)=-\frac{5}{2}h'(0);~
 I_{\rm {Gr,b}}=-5h'(0){\rm Vol}_{\partial M}.
\end{equation}
 Then we obtain
 \begin{eqnarray}
\widetilde{{\rm Wres}}[(\pi^+D^{-2}_{T})^2]_{i}&=&0;\\
\widetilde{{\rm Wres}}[(\pi^+D^{-2}_{T})^2]_{b}&=&\int_{\partial M}\Phi=Q_{0}\pi\Omega_{5}{\rm Vol}_{\partial M}.\nonumber\\
\end{eqnarray}
where\begin{eqnarray}
Q_{0}&=&\frac{-1}{2}s_{M}(x_{0})+\frac{35}{24}s_{\partial_{M}}(x_{0})-2T^{2}|V|^{2}-(\frac{45}{2}+\frac{5i}{4})h'(0)\nonumber\\
&&+\big(\frac{45i}{32}-\frac{3947}{384}\big)\big(h'(0)\big)^{2}
+\big(\frac{247}{32}-\frac{45i}{32}\big)h''(0)
\end{eqnarray}
By (4.5)-(4.6), we obtain
\begin{cor}
Let $M$ be a $7$-dimensional
compact spin manifold with the boundary $\partial M$ and the metric
$g^M$ as above and  $D_{T}$ be the Witten deformation on $\widehat{M}$, then
\begin{eqnarray}
I_{\rm {Gr,b}}=-\frac{5h'(0)}{Q_{0}\pi \Omega_{5}}\widetilde{{\rm Wres}}[(\pi^+D^{-1}_{T})^2]_{b}.
\end{eqnarray}
\end{cor}

\section*{ Acknowledgements}
On behalf of all authors, the corresponding author states that there is no conflict of interest.
This work is supported by NSFC(11901322) and Inner Mongolia Natural Science Foundation(2018LHO1004).


\section*{References}
\begin{thebibliography}{00}

\bibitem{Gu} V. W. Guillemin.: A new proof of Weyl's formula on the asymptotic distribution of eigenvalues. Adv. Math. 55, no. 2, 131-160, (1985).
\bibitem{Wo} M.Wodzicki.: local invariants of spectral asymmetry. Invent. Math. 75(1), 143-178, (1995).
\bibitem{Co1} A. Connes.: Quantized calculus and applications.  XIth International Congress of Mathematical Physics(Paris,1994),
Internat Press, Cambridge, MA, 15-36, (1995).
\bibitem{Co2} A. Connes.: The action functinal in Noncommutative geometry. Comm. Math. Phys. 117, 673-683, (1998).
\bibitem{Ka} D. Kastler.: The Dirac Operator and Gravitation. Comm. Math. Phys. 166, 633-643, (1995).

\bibitem{KW} W. Kalau and M. Walze.: Gravity, Noncommutative geometry and the Wodzicki residue. J. Geom. Physics. 16, 327-344,(1995).
\bibitem{Ac} T. Ackermann.: A note on the Wodzicki residue. J. Geom.Phys. 20, 404-406, (1996).
\bibitem{Po} R. Ponge.: Noncommutative geometry and lower dimensional volumes in Riemannian geometry. Letters in Mathematical Physics.
83, no.1, 19-32, (2008).
\bibitem{FGLS} B. V. Fedosov, F. Golse, E. Leichtnam, E. Schrohe.: The noncommutative residue for manifolds with boundary. J. Funct. Anal.
142, 1-31, (1996).

\bibitem{Wa1} Y. Wang.: Diffential forms and the Wodzicki residue for Manifolds with Boundary. J. Geom. Physics. 56, 731-753, (2006).
\bibitem{Wa2} Y. Wang.: Diffential forms the Noncommutative Residue for Manifolds with Boundary in the non-product Case,
      Letters in Mathematical Physics. 77, 41-51, (2006).
\bibitem{Wa3} Y. Wang.: Gravity and the Noncommutative Residue for Manifolds with Boundary. Letters in Mathematical Physics. 80, 37-56, (2007).
\bibitem{WJ2} J. Wang, Y. Wang, C. L. Yang. A Kastler-Kalau-Walze Type Theorem for 7-Dimensional Manifolds with Boundary. Abstract and Applied Analysis. 2014, Art. ID 465782, 1-18 (2014).
\bibitem{Wa4} Y. Wang.: Lower-Dimensional Volumes and Kastler-kalau-Walze Type Theorem for Manifolds with Boundary .
\bibitem{WJ1} J. Wang, Y. Wang, C. L. Yang. : Dirac operators with torsion and the noncommutative residue for manifolds with boundary.
     J. Geom. Physics. 81, 92-111, (2014).
\bibitem{wpz} W. Zhang. Lectures on chern-weil theory and witten deformation, vol.4. World  Scientific  Publishing  Co. Pte. Ltd., 2001.
      Commun. Theor. Phys. Vol 54, 38-42, (2010).
\bibitem{Y} Y. Yu.: The Index Theorem and The Heat Equation Method, Nankai Tracts in Mathematics-Vol.2, World Scientific Publishing, (2001).
\bibitem{BGV} Y. Yu.:N. Berline; E. Getzler; M. Vergne, Heat kernels and Dirac operators. Springer- Verlag, Berlin, 1992.
\bibitem{U}  W. J. Ugalde, Differential forms and the Wodzicki residue, arXiv: Math, DG/0211361.
\bibitem{KH1}Kai Hua Bao, Ying Lei. The Kastler¨CKalau¨CWalze type theorem about Witten
deformation for manifolds with boundary. Journal of Geometry and Physics, May 2018, Volume 131(2018). 170-181.
\bibitem{KH2}Kai Hua BAO, Ai Hui SUN, Chao DENG. The Noncommutative Residue about Witten Deformation[J]. Acta Mathematica Sinica, English Series, 2019,35(4): 550-568.
\bibitem{KH3}Kai Hua Bao,Ai Hui Sun, Jian Wang. A Kastler-Kalau-Walze Type Theorem for 7-dimensional Spin Manifolds with Boundary about Dirac Operators with Torsion. Journal of Geometry and Physics,December 2016, Volume 110. 213¨C232.


%%%%%%%%%%%%%%%%%%%%%%%%%%%%%%%%%%%%%%%%%%%%%%%%%%%%%%%%%%%%%%%%%%%%%%%%%%%%%%%%%%%%%%%%%%%%%%%%%%%%%%%%%%%%%%%%%%%%%%%%%%%%%%%%%%%%%%%%%%%%

\bibitem{MA} M. Adler.: On a trace functional for formal pseudo-differential operators and the symplectic
structure of Korteweg-de Vries type equations, Invent. Math. 50, 219-248,(1979).
\bibitem{S} E. Schrohe.:  Noncommutative residue, Dixmier's trace, and heat trace expansions on manifolds with boundary. Contemp. Math.
242, 161-186, (1999).
\bibitem{TJ} T. Ackermann, J. Tolksdorf, A generalized Lichnerowicz formula, the Wodzicki residue and gravity, J. Geom. Phys. 19 (1996) 143-150.
 \bibitem{FC1} F. Pf\"{a}ffle, C.A. Stephan, On gravity, torsion and the spectral action principle, J. Funct. Anal. 262 (2012) 1529-1565.
\bibitem{FC2} F. Pf\"{a}ffle, C.A. Stephan, Chiral asymmetry and the spectral action, Comm. Math. Phys. 321 (2013) 283-310.

\end{thebibliography}
\end{document}